\documentclass[reqno,draft]{amsart}

\usepackage{amssymb}
\usepackage{pst-all}

\renewcommand{\d}[1]{\,\mathrm{d}#1}
\newcommand{\R}{\mathbb{R}}
\newcommand{\K}{\mathbb{K}}
\newcommand{\C}{\mathbb{C}}
\newcommand{\N}{\mathbb{N}}
\newcommand{\Z}{\mathbb{Z}}
\newcommand{\Int}{\mathop\mathrm{Int}\nolimits}

\newcommand{\ind}{\mathop\mathrm{ind}\nolimits}

\newcommand{\Ker}{\mathop\mathrm{Ker}\nolimits}
\newcommand{\coker}{\mathop\mathrm{coKer}\nolimits}
\renewcommand{\Im}{\mathop\mathrm{Im}\nolimits}
\newcommand{\sign}{\mathop\mathrm{sign}\nolimits}
\newcommand{\dist}{\mathop\mathrm{dist}\nolimits}

\newcommand{\B}{\mathcal{B}}

\newcommand{\bde}{\begin{definition}\label}
\newcommand{\ede}{\end{definition}}
\newcommand{\bpr}{\begin{proposition}\label}
\newcommand{\epr}{\end{proposition}}
\newcommand{\ble}{\begin{lemma}\label}
\newcommand{\ele}{\end{lemma}}
\newcommand{\bco}{\begin{corollary}\label}
\newcommand{\eco}{\end{corollary}}
\newcommand{\bre}{\begin{remark}\label}
\newcommand{\ere}{\end{remark}}
\newcommand{\beq}{\begin{equation}\label}
\newcommand{\eeq}{\end{equation}}
\newcommand{\ben}{\begin{enumerate}}
\newcommand{\een}{\end{enumerate}}
\newcommand{\bit}{\begin{itemize}}
\newcommand{\eit}{\end{itemize}}

\newcommand{\bto}{\left(\begin{array}}
\newcommand{\eto}{\end{array}\right)}

\newcommand{\bgr}{\left\{\begin{array}}
\newcommand{\enu}{\end{array}\right.}

\newcommand{\bnu}{\left.\begin{array}}

\theoremstyle{plain}
\newtheorem{theorem}{Theorem}[section]
\newtheorem{corollary}[theorem]{Corollary}
\newtheorem{lemma}[theorem]{Lemma}
\newtheorem{proposition}[theorem]{Proposition}
\theoremstyle{definition}
\newtheorem{definition}[theorem]{Definition}
\newtheorem{remark}[theorem]{Remark}
\newtheorem{example}[theorem]{Example}

\numberwithin{equation}{section}

\newcommand{\cc}{\overline{\mbox{\textup{co}}}\,}

\newcommand{\eps}{\varepsilon}

\begin{document}

\title{A New Spectrum for Nonlinear Operators in Banach Spaces}

\author[A.\ Calamai]{Alessandro Calamai}
\author[M.\ Furi]{Massimo Furi}
\author[A.\ Vignoli]{Alfonso Vignoli}

\date{}

\address{
Alessandro Calamai,
Dipartimento di Scienze Matematiche, Universit\`a Politecnica delle
  Marche, Via Brec\-ce Bianche,
I-60131 Ancona, Italy.
E-mail addresses: calamai@math.unifi.it, calamai@dipmat.univpm.it}
\address{
Massimo Furi,
Dipartimento di Matematica Applicata ``G.~Sansone'', Via S.\ Mar\-ta 3,
I-50139 Firenze, Italy.
E-mail address: massimo.furi@unifi.it}
\address{
Alfonso Vignoli,
Dipartimento di Matematica, Universit\`a di Roma ``Tor Vergata'', Via
della Ricerca Scientifica,
I-00133 Roma, Italy.
E-mail address: vignoli@mat.uniroma2.it}

\thanks{We thank both the referees for insightful comments that helped
us to improve the paper significantly}

\begin{abstract}
Given any continuous self-map $f$ of a Banach space $E$ over $\K$
(where $\K$ is $\R$ or $\C$) and given any point $p$ of $E$, we define
a subset $\sigma(f,p)$ of $\K$, called \emph{spectrum of $f$ at $p$},
which coincides with the usual spectrum $\sigma(f)$ of $f$ in the
linear case.  More generally, we show that $\sigma(f,p)$ is always
closed and, when $f$ is $C\sp 1$, coincides with the spectrum
$\sigma(f'(p))$ of the Fr\'echet derivative of $f$ at $p$. Some
applications to bifurcation theory are given and some peculiar
examples of spectra are provided.
\end{abstract}

\maketitle


\section{Introduction}
\label{Introduction}
\setcounter{equation}{0}

In spite of the fact that this paper deals with spectral theory for
merely continuous nonlinear operators, we feel that the special case
of (Fr\'echet) differentiable operators is worth a particular regard.

In Nonlinear Analysis the class of differentiable
operators is very natural and important. It is therefore not
surprising that one of the first attempts of defining a spectrum for
nonlinear operators was undertaken for this type of maps. Namely,
already in 1969, J.W.\ Neuberger proposed a simple definition of
spectrum that runs as follows (see~\cite{Neu}).
Given a Banach space $E$ over $\K$
(where $\K$ is $\R$ or $\C$)
and a map $f:E \to E$ of class $C\sp 1$, the set
\[
\rho_N(f)=\{\lambda\in \K : \lambda I - f
\mbox{ is bijective and } (\lambda I-f)\sp{-1}
\mbox{ of class } C\sp 1 \}
\]
is called the \textit{Neuberger resolvent set} and its complement
\[
\sigma_N(f)= \K \backslash \rho_N(f)
\]
the \textit{Neuberger spectrum} of $f$.
As pointed out in~\cite{ApDo}, this definition is interesting since it
agrees with the usual one in the linear case.
Moreover, the Neuberger spectrum is always
nonempty if the underlying space is complex.
However, it has some drawbacks. In particular, the Neuberger spectrum
need not be closed.
This already happens in the one dimensional case
(see e.g.~\cite{ADPV}).

Our approach here is completely different.
Given an open subset $U$ of $E$, a continuous map $f:U \to E$,
and a point $p\in U$, we introduce
the concept of \textit{spectrum of the map $f$ at the point $p$},
which we denote by $\sigma(f,p)$.
Our spectrum is close in spirit to the asymptotic spectrum introduced
in~\cite{fumavi78} for continuous maps defined on the whole space $E$.
This asymptotic spectrum turns out to be closed,
and coincides with the usual spectrum in the linear case.
While the asymptotic spectrum is related to the asymptotic behavior
of a map, our spectrum $\sigma(f,p)$
depends only on the germ of $f$ at $p$.

The first attempt to introduce a notion of spectrum at a point was
undertaken in~\cite{May} by means of a suitable local adaptation of
Neuberger's ideas.
Subsequently, the last two authors gave
in~\cite{fuvi77} another definition of spectrum at a point.
Namely, given $f$ and $p$ as above, they defined in~\cite{fuvi77} a
spectrum $\Sigma(f,p)$
which, in the case of a bounded linear operator
$L: E\to E$, gives only a part of the classical spectrum $\sigma(L)$.
Namely, $\Sigma(L,p)$ reduces to the approximate point spectrum of
$L$.
Therefore, the definition of $\Sigma(f,p)$ is somehow nonexhaustive.
For example, if $L : \ell \sp 2 (\C) \to \ell \sp 2 (\C)$ is the
\textit{right-shift operator}, $\sigma(L)$ is the unit disk of $\C$
and the approximate point spectrum of $L$ coincides with
$\partial \sigma(L) = S \sp 1$.

We are now in the position to fill the gap.
For this purpose we use the topological concept of
\textit{zero-epi map} (see~\cite{fumavi80})
as well as two numerical characteristics recently
introduced in~\cite{Ca} (see also~\cite{CaPhd}).
These are the notions which turn out to be crucial in our new
definition of spectrum at a point.

Given $f$ and $p$ as above, we show that $\sigma(f,p)$ is always
closed and contains $\Sigma(f,p)$.
Moreover, for a bounded linear operator $L$, we get
$\sigma(L,p) = \sigma(L)$ for any $p \in E$,
where $\sigma(L)$ is the usual spectrum of $L$.
More precisely, if a map $f$ is $C\sp 1$, then $\sigma(f,p)$ is equal to
the spectrum of the Fr\'echet derivative $\sigma(f'(p))$.
Thus, for $C\sp 1$ maps and $\K =\C$, our spectrum is always nonempty as
in the Neuberger case.

As we shall see, our spectrum at a point turns out to be useful to
tackle bifurcation problems in the non-differentiable case.

Our paper closes with a few examples illustrating some peculiarities
of the nonlinear case.


\section{Notation and preliminary results}
\label{sec-mnc}

In this section we collect some notions that we will need in order to
define our spectrum at a point.

\medskip

Throughout the paper, $E$ and $F$ will be two Banach spaces over $\K$,
where $\K$ is $\R$ or $\C$.
Given a subset $X$ of $E$, by $C(X,F)$ we denote the set of all
continuous maps from $X$ into $F$.

Let $U$ be an open subset of $E$,
$f\in C(U,F)$ and $p\in U$.
Denote by $U_p$ the open neighborhood $\{x\in E: p+x\in U\}$ of
$0\in E$, and define $f_p\in C(U_p,F)$ by
\[
f_p(x)=f(p+x)-f(p).
\]

Given a subset $A$ of a metric space $X$, we will denote,
respectively, by $\overline A$, $\Int A$ and $\partial A$ the
closure, the interior and the boundary of $A$.

Let $f:X\to Y$ be a continuous map between two metric spaces.
We recall that $f$ is said to be \textit{compact} if $f(X)$ is
relatively compact, and \textit{completely continuous} if it is compact
on any bounded subset of $X$. If for any $p\in X$ there exists a
neighborhood $V$ of $p$ such that the restriction $f|_V$ is compact,
then $f$ is called \textit{locally compact}.
The map $f$ is said to be \textit{proper} if $f\sp{-1}(K)$ is compact for
any compact subset $K$ of $Y$.
Recall that a proper map sends closed sets into closed sets.
Given $p\in X$, the map
$f$ is \textit{locally proper at $p$} if there exists
a closed neighborhood $V$ of $p$ such that the restriction $f|_V$ is
proper.

In the sequel we will adopt the conventions
$\sup \emptyset = -\infty$ and $\inf \emptyset = + \infty$,
so that if $A\subseteq B\subseteq \R$, then $\inf A\geq \inf B$ and
$\sup A\leq \sup B$ even when $A$ is empty.


\subsection{Definition and properties of $\alpha_p(f)$ and
   $\omega_p(f)$}

Let us recall the definition and properties of the
Kuratowski measure of noncompactness (see~\cite{Ku}).

The \emph{Kuratowski measure of noncompactness} $\alpha(A)$ of a subset $A$ of
$E$ is defined as the infimum of real numbers $d>0$ such that $A$ admits a
finite covering by sets of diameter less than $d$. In particular, if $A$ is
unbounded, we have $\alpha(A) = \inf \emptyset = +\infty$.
Notice that, if $E$ is finite dimensional, then $\alpha(A) = 0$ for any
bounded subset $A$ of $E$.

We summarize the following properties of the measure of
noncompactness. Given a subset $A$ of $E$, we denote by $\cc A$
the closed convex hull of $A$, and by $[0,1] A$ the set
$\{ \lambda x : \lambda\in [0,1], x\in A\}$.

\begin{proposition}
\label{prop-mnc}
Let $A$ and $B$ be subsets of $E$. Then
\begin{enumerate}
\item $\alpha(A) = 0$ if and only if $\overline A$ is
    compact; \smallskip
\item $\alpha(\lambda A) = |\lambda| \alpha(A)$ for any
    $\lambda \in \R$; \smallskip
\item $\alpha(A+B) \leq \alpha(A) + \alpha(B)$; \smallskip
\item if $A \subseteq B$, then $\alpha(A) \leq
    \alpha(B)$; \smallskip
\item $\alpha(A \cup B) =\max \{ \alpha(A) , \alpha(B) \}$;
    \smallskip
\item $\alpha([0,1] A) =\alpha(A)$; \smallskip
\item $\alpha(\cc A ) =\alpha(A)$.
\end{enumerate}
\end{proposition}

Properties~(1)--(6) are straightforward consequences of the
definition, while the last one is due to Darbo~\cite{Dar}.

\medskip

Given a subset $X$ of $E$ and $f\in C(X,F)$, we recall
the definition of the following two extended real numbers (see
e.g.~\cite{fumavi78}) associated with the map $f$:
\[
\alpha(f) = \sup \left\{ \frac{\alpha (f(A))}{\alpha (A)} : \,
A\subseteq X \mbox{ bounded, } \alpha (A) > 0 \right\},
\]
and
\[
\omega(f) = \inf \left\{ \frac{\alpha (f(A))}{\alpha (A)} : \,
A\subseteq X \mbox{ bounded, } \alpha (A) > 0 \right\},
\]
where $\alpha$ is the Kuratowski measure of noncompactness
(in~\cite{fumavi78} $\omega(f)$ is denoted by $\beta(f)$, however, we
prefer here the more recent notation $\omega(f)$ as in~\cite{EW}).
It is convenient to put $\alpha(f)=-\infty$ and $\omega(f)=+\infty$
whenever $E$ is finite dimensional and $f\in C(X,F)$.
This agrees with the notation
$\sup \emptyset = -\infty$ and $\inf \emptyset = + \infty$,
and allows to simplify some statements.

It is important to observe that $\alpha(f) \leq 0$ if and only if $f$
is completely continuous.
Moreover, $\omega(f) >0$ only if $f$  is proper on bounded closed
sets.
For a comprehensive list of properties of $\alpha(f)$ and $\omega(f)$
we refer to~\cite{fumavi78}.

\medskip

Let $U$ be an open subset of $E$, $f\in C(U,F)$ and $p\in U$.
We recall the definitions of $\alpha_p(f)$ and $\omega_p(f)$ given
in~\cite{Ca}. Roughly speaking, these numbers are the local analogues
of $\alpha(f)$ and $\omega(f)$.

Let $B(p,r)$ denote the open ball in $E$ centered at $p$ with radius
$r>0$. Suppose that $B(p,r)\subseteq U$ and consider the number
\[
\alpha(f|_{B(p,r)})= \sup \left\{ \frac{\alpha(f(A))}{\alpha(A)} :
A \subseteq B(p,r),\; \alpha(A) > 0 \right\},
\]
which is nondecreasing as a function of $r$. Hence, we can define
\[
\alpha_p(f)=\lim_{r\to 0} \alpha(f|_{B(p,r)}).
\]
Clearly, $\alpha_p(f) \leq \alpha(f)$.
Analogously, define
\[
\omega_p(f)=\lim_{r\to 0} \omega (f|_{B(p,r)}).
\]
Obviously, $\omega_p(f) \geq \omega(f)$.
Notice that $\alpha_p(f)=\alpha_0(f_p)$ and
$\omega_p(f)=\omega_0(f_p)$.

Let us point out that, if $E$ is finite dimensional, then
$\alpha_p(f)=-\infty$ and $\omega_p(f)=+\infty$ for any $p\in U$.

\medskip

With only minor changes, it is easy to show that the main properties of
$\alpha$ and $\omega$ hold for $\alpha_p$ and $\omega_p$ as well.
In fact, the following propositions hold (see~\cite{Ca}).

\begin{proposition} \label{prop-alfaomega}
Let $E$ be infinite dimensional.
Given an open subset $U$ of $E$, $f,g\in C(U,F)$ and
$p\in U$, one has
\begin{enumerate}
\item $\alpha_p (c f) = |c| \alpha_p (f)$ and
$\omega_p (c f) = |c| \omega_p (f)$, for any $c \in
\K$;\smallskip
\item $\omega_p (f) \leq \alpha_p (f)$; \smallskip
\item $|\alpha_p (f) - \alpha_p (g)| \leq \alpha_p (f+g) \leq
   \alpha_p(f)+\alpha_p (g)$; \smallskip
\item $\omega_p (f) - \alpha_p (g) \leq \omega_p (f+g) \leq
   \omega_p(f) + \alpha_p (g)$; \smallskip
\item if $f$ is locally compact, then $ \alpha_p (f) =0 $; \smallskip
\item if $\omega_p (f) >0$, then $f$ is locally proper at $p$;
   \smallskip
\item if $f$ is a local homeomorphism and $\omega_p (f)>0$, then
$\alpha_q (f\sp{-1}) \omega_p (f)=1$, where $q=f(p)$.
\end{enumerate}
\end{proposition}

\begin{proposition} \label{prop-alfacomposto}
Let $E$, $F$ and $G$ be infinite dimensional Banach spaces,
$U \subseteq E$ and $V \subseteq F$ open,
$f \in C(U,F)$ and $g \in C(V,G)$.
Given $p \in f\sp{-1}(V)$, let $q = f(p) \in V$. Then, for the composite map
$gf$, we have
\begin{enumerate}
\item $\alpha_p(gf) \leq \alpha_q(g) \alpha_p(f)$; \smallskip
\item $\omega_q(g) \omega_p(f) \leq \omega_p(gf) \leq \alpha_q(g)\omega_p(f)$.
\end{enumerate}
\end{proposition}

\begin{proposition}
\label{prop-alfaposom}
If $f:E\to F$ is positively homogeneous, then
$\alpha_0(f)=\alpha(f)$ and $\omega_0(f)=\omega(f)$.
\end{proposition}

Clearly, for a bounded linear operator $L: E\to F$, the numbers
$\alpha_p(L)$ and $\omega_p(L)$ do not depend on the point $p$ and
coincide, respectively, with $\alpha(L)$ and $\omega(L)$.
Furthermore, for the $C\sp 1$ case the following result holds.

\begin{proposition}[\cite{Ca}]
\label{prop-lin2}
Let $f:U\to F$ be of class $C\sp 1$. Then, for any $p\in U$ we
have $\alpha_p(f)=\alpha(f'(p))$ and $\omega_p(f)=\omega(f'(p))$.
\end{proposition}


\subsection{Definition and properties of $|f|_p$ and $d_p(f)$}

As before, let $f \in C(U,F)$, and fix $p\in U$.
Following~\cite{fuvi77}, define
\[
|f|_p = \limsup_{x\to 0} \frac{\|f(p+x)-f(p)\|}{\|x\|}
\]
and
\[
d_p(f) = \liminf_{x\to 0} \frac{\|f(p+x)-f(p)\|}{\|x\|}.
\]
Notice that $|f|_p=|f_p|_0$ and $d_p(f)=d_0(f_p)$.

Following~\cite{fuvi77}, the map $f$ will be called
\textit{quasibounded at $p$} if $|f|_p < +\infty$.

\medskip

The next proposition collects some useful properties of $|f|_p$
and $d_p(f)$.

\begin{proposition}[\cite{fuvi77}]
\label{prop-dq-proprieta}
Given an open subset $U$ of $E$, $f,g\in C(U,F)$ and
$p\in U$, one has
\begin{enumerate}
\item $|c f|_p = |c| |f|_p$ and
$d_p (c f) = |c| d_p (f)$, for any $c \in
\K$;\smallskip
\item $0\leq d_p (f) \leq |f|_p$; \smallskip
\item $|d_p (f) - d_p (g)| \leq |f-g|_p$; \smallskip
\item $|f+g|_p \leq |f|_p+|g|_p$; \smallskip
\item $d_p (f) - |g|_p \leq d_p (f+g) \leq d_p(f) + |g|_p$.
\end{enumerate}
\end{proposition}

\begin{remark}
\label{oss-dq-posom}
If $f:E\to F$ is positively homogeneous, then
$|f|_0$ coincides with the quasinorm of $f$, i.e.\
\[
|f|=\limsup_{\|x\|\to +\infty} \frac{\|f(x)\|}{\|x\|},
\]
and $d_0(f)$ coincides with the number
\[
d(f) = \liminf_{\|x\|\to +\infty} \frac{\|f(x)\|}{\|x\|}.
\]
More precisely, one has
\[
|f|_0=|f|=\sup_{\|x\|=1} \|f(x)\|
\quad \mbox{ and } \quad
d_0(f)=d(f)=\inf_{\|x\|=1} \|f(x)\|.
\]
\end{remark}

For the $C\sp 1$ case the following result holds.

\begin{proposition}[\cite{fuvi77}] \label{prop-lin2-dq}
Let $f:U\to F$ be of class $C\sp 1$. Then, for any $p\in U$ we
have $|f|_p=\|f'(p)\|$ and $d_p(f)=d(f'(p))$.
\end{proposition}

\medskip

The next result was proved in~\cite{fumavi78}
(see Proposition~3.2.3). Here we give a different proof.

\begin{proposition}
\label{prop-betaposom}
Let $f:E\to F$ be positively homogeneous, and assume that
$\omega(f)>0$.
Then, the equation $f(x)=0$ has only the trivial solution if and only
if $d(f)>0$.
\end{proposition}

\begin{proof}
Clearly, $d(f)>0$ implies that the equation
$f(x)=0$ has only the trivial solution.
Thus, it is enough to show that the
conditions $\omega(f)>0$ and $d(f)=0$ imply the existence of an
element $x \in E$, $x \neq 0$, such that $f(x)=0$.
Now, $\omega(f)>0$ implies that $f$ is proper on the set
$A= \{ x \in E : \|x\|=1 \}$.
Therefore, $f(A)$ is closed.
Hence, from $d(f)=0$ it follows that $\dist(0,f(A))=0$.
Consequently, there exists $x \in A$ such that $f(x)=0$
and the assertion follows.
\end{proof}

Evidently, for a bounded linear operator $L: E\to F$, the number
$|L|_p$ does not depend on the point $p$ and coincides with the norm
$\|L\|$.
Analogously, $d_p(L)$ is independent of $p$ and coincides with
$d(L)$.

\medskip

The following proposition (see e.g.~\cite{fumavi78})
gives information on the numbers $d$,
$\alpha$ and $\omega$ in the context of linear operators.
Given a bounded linear operator $L: E \to F$, by $L\sp{*}$ we denote the
adjoint of $L$.

\begin{proposition} \label{prop-lin1-dq}
Let $L:E\to F$ be a bounded linear operator.
We have
\begin{enumerate}
\item $\alpha(L)\leq\|L\|$;\smallskip
\item $d(L)\leq\omega(L)$;\smallskip
\item $d(L)>0$ if and only if
$L$ is injective and its inverse is continuous; \smallskip
\item $d(L)\sp{-1}=\|L\sp{-1}\|$ for any linear isomorphism; \smallskip
\item $\omega(L)>0$ if and only if $L$ is left semi-Fredholm; that is,
$\Im L$ is closed and $\dim \Ker L < + \infty$; \smallskip
\item $\omega(L\sp*)>0$ if and only if $L$ is right semi-Fredholm; that is,
$\Im L$ is closed and $\dim \coker L < + \infty$.
\end{enumerate}
\end{proposition}

As a consequence of the previous proposition, $L$ is a
\textit{Fredholm operator}
(that is, $\Ker L$ and $\coker L$ are finite dimensional and $\Im L$
is closed)
if and only if $\omega(L)>0$ and $\omega(L\sp*)>0$.


\section{zero-epi maps at a point}

In this section we recall the definition and properties of zero-epi
(and $y$-epi)
maps and we introduce a local analogue of this concept.


\subsection{Zero-epi maps}

Let $U$ be a bounded open subset of $E$, and
$f\in C(\overline U,F)$.
We recall the following definitions given in~\cite{fumavi80}.

\begin{definition}
Given $y\in F$, we say that $f$ is \textit{y-admissible (on $U$)} if
$f(x)\neq y$ for any $x\in \partial U$.
\end{definition}

\begin{definition}
We say that $f$ is \textit{$y$-epi (on $U$)} if it is $y$-admissible
and for any compact map $h:\overline U \to F$ such that $h(x)=y$ for
all $x \in \partial U$ the equation $f(x)=h(x)$ has a solution in
$U$.
\end{definition}

Notice that $f$ is $y$-epi if and only if the map $f-y$, defined by
$(f-y)(x)=f(x)-y$, is $0$-epi (zero-epi).

\medskip

The main properties of zero-epi maps are analogous to some
of the properties which characterize the Leray--Schauder degree.
In particular, we mention the following one to be used in the sequel.

\begin{theorem}[Homotopy invariance]
\label{omo-epi}
Let $f\in C(\overline U,F)$ be $0$-epi on $U$, and
$H:\overline U \times [0,1] \to F$ be compact and such that $H(x,0)=0$
for any $x \in \overline U$. Assume that $f(x)+H(x,\lambda)\neq 0$ for
any $x \in \partial U$ and any $\lambda \in [0,1]$.
Then, the map $f(\cdot)+H(\cdot,1)$ is $0$-epi on $U$.
\end{theorem}

As a straightforward consequence of the above homotopy invariance
property, we obtain the following Rouch\'e--type perturbation result for
zero-epi maps.

\begin{corollary}
\label{omo-epi-cor}
Let $f \in C(\overline U,F)$ be $0$-epi on $U$, and
$k \in C(\overline U,F)$ be compact and such that
$\|k(x)\| < \|f(x)\|$ for any $x \in \partial U$.
Then, the equation $f(x)=k(x)$ has a solution in $U$.
\end{corollary}

\begin{proof}
Apply Theorem~\ref{omo-epi} with $H(x,\lambda)= - \lambda k(x)$.
\end{proof}

The next local surjectivity property of zero-epi maps has been
proved in~\cite{fumavi80}.

\begin{theorem}[Local surjectivity]
\label{epi-onto}
Let $f\in C(\overline U,F)$ be $0$-epi on $U$ and proper on
$\overline U$. Then, $f$ maps $U$ onto a neighborhood of the origin.
More precisely, if $V$ is the connected component of
$F \backslash f(\partial U)$ containing the origin, then
$V \subseteq f(U)$.
\end{theorem}

\medskip

The next properties of zero-epi maps are straightforward consequences
of the definition.

\begin{proposition}[Localization] \label{prop-local1}
If $f\in C(\overline U,F)$ is $0$-epi on $U$, and $U_1$ is an open
subset of $U$ containing $f\sp{-1}(y)$, then $f|_{\overline U_1}$ is
$0$-epi.
\end{proposition}

\begin{proposition}[Topological invariance] \label{prop-topol}
Let $U$ and $V$ be bounded open subsets of $E$, and
$\varphi : \overline U \to \overline V$ be a homeomorphism
such that $\varphi(\partial U) = \partial V$.
Given $f\in C(\overline U,F)$, define $g\in C(\overline V,F)$ by
$g=f\varphi\sp{-1}$.
Then, $f$ is $0$-admissible on $U$ if and only if
$g$ is $0$-admissible on $V$,
and $f$ is $0$-epi on $U$ if and only if
$g$ is $0$-epi on $V$.
\end{proposition}

\begin{proposition}
Let $f\in C(\overline U,F)$ be $0$-epi on $U$, and
$c \in \K \backslash \{0\}$.
Then, the map $cf$ is $0$-epi on $U$.
\end{proposition}

The following proposition shows that if $f\in C(\overline U,E)$ is a
compact vector field (i.e., a compact perturbation of the identity)
and $U$ is a ball, the converse of the localization
property stated in Proposition~\ref{prop-local1} holds.
This converse property is known to hold in the more general case in
which $U$ is a bounded open set and $f = I - g$, where $g$ is an
$\alpha$-contraction (see~\cite{GV}).
Our assumptions on $U$ and $f$ allow us to provide an independent and
elementary proof, which is not based on the Hopf classification theorem
in degree theory.

By $B_r$ we mean the open ball in $E$ about the origin with radius
$r>0$.

\begin{proposition}[Antilocalization]
\label{prop-taglio}
Let $U_1 = B_{r_1}$ and $g\in C(\overline U_1,E)$ be compact.
Assume that there exists $0<r_0<r_1$ such that $I-g$ is $0$-epi on
$U_0=B_{r_0}$, and that the equation $x-g(x)=0$ has no solutions
in $\overline U_1 \backslash U_0$.
Then, $I-g$ is $0$-epi on $U_1$.
\end{proposition}

\begin{proof}
First, we claim that the $0$-admissible map
\[
x \mapsto x - \frac{r_1}{r_0} g \left( \frac{r_0}{r_1} x \right)
\]
is $0$-epi on $U_1$.
To see this, let $h:\overline U_1 \to E$
be compact and such that $h(x)=0$ for
all $x \in \partial U_1$.
We have to prove that the equation
\begin{equation} \label{equ-taglio1}
x-\frac{r_1}{r_0} g \left( \frac{r_0}{r_1} x \right)=h(x)
\end{equation}
has a solution in $U_1$.
This equation is equivalent to
\[
\frac{r_0}{r_1} x - g \left( \frac{r_0}{r_1} x \right)=
\frac{r_0}{r_1} h(x), \quad x \in U_1.
\]
Put $y=\frac{r_0}{r_1}x$.
Since $x \in U_1$ if and only if
$y \in U_0$, we get
\begin{equation} \label{equ-taglio2}
y - g (y)= \frac{r_0}{r_1} h \left( \frac{r_1}{r_0} y \right),
\quad y \in U_0.
\end{equation}
Notice that the map $y \mapsto h \left( \frac{r_1}{r_0} y \right)$ is
compact and $h \left( \frac{r_1}{r_0} y \right) =0$ for
all $y \in \partial U_0$.
Since $I-g$ is $0$-epi on $U_0$,
equation~\eqref{equ-taglio2} has a solution $\bar y \in U_0$.
Consequently, $\bar x = \frac{r_1}{r_0} \bar y \in U_1$
is a solution of equation~\eqref{equ-taglio1}.
Thus, the map
$x \mapsto x - \frac{r_1}{r_0} g \left( \frac{r_0}{r_1} x \right)$
is $0$-epi on $U_1$, as claimed.

Now, consider the map
\[
H:\overline U_1 \times \left[ \frac{r_0}{r_1},1 \right] \to E, \quad
H(x,s)= - \frac 1s g(sx).
\]
Clearly, $H$ is compact and
$x \mapsto x + H(x,\frac{r_0}{r_1})$ is $0$-epi on $U_1$.
Since $H(\cdot,1) = -g$,
by the homotopy invariance property of zero-epi maps,
the assertion
follows if we show that $x+H(x,s)\neq 0$ for
any $x \in \partial U_1$ and any $s \in [\frac{r_0}{r_1},1]$.
To prove this, consider the equation
\[
x - \frac 1s g(sx) = 0, \quad x \in \overline U_1, \
s \in \left[ \frac{r_0}{r_1},1 \right],
\]
which is equivalent to
\begin{equation} \label{equ-taglio3}
sx - g(sx) = 0, \quad x \in \overline U_1, \
s \in \left[ \frac{r_0}{r_1},1 \right].
\end{equation}
Observe that from $\|x\| = r_1$ it follows that
$r_0 \leq \|sx\| \leq r_1$.
Thus, by assumption, the equation~\eqref{equ-taglio3}
has no solution $x \in \partial U_1$ for all
$\frac{r_1}{r_0} \leq s \leq 1$.
Consequently, by the homotopy invariance property of
zero-epi maps, the map $I-g$ is $0$-epi on $U_1$.
This completes the proof.
\end{proof}

\medskip

The following coincidence theorem for zero-epi maps is due to V\"ath
(see~\cite{Vaeth02}).
Here we give an independent proof, following closely some ideas
contained in Theorem~4.2.1 of~\cite{fumavi78}.

\begin{theorem}[Coincidence theorem for zero-epi maps]
\label{teo-perturba}
Let $U$ be a bounded open subset of $E$, $f\in C(\overline U,F)$
a $0$-epi map with $\omega(f)>0$, and $k\in C(\overline U,F)$.
Assume that $\alpha(k)<\omega(f)$ and
\[
\sup \{\|k(x)\|: x\in \overline U \}
< \inf \{\|f(x)\|: x\in \partial U \}.
\]
Then the equation
$f(x)=k(x)$ has a solution in $U$.
\end{theorem}

\begin{proof}
Define, by induction, the following sequence of closed subsets of $E$:
\[
X_0=\overline U, \qquad X_{n+1}=f\sp{-1}(\cc(k(X_n))), \ n\geq 0.
\]
Let us show first that $X_1\neq\emptyset$.
To this purpose, notice that the image $f(U)$ contains the open ball
$B_r$ in $F$ centered at the origin with radius
\[
r = \inf \{\|f(x)\|: x\in \partial U \}>0.
\]
Indeed, let $y\in F$ be such that $\|y\|< r$.
Then, from Corollary~\ref{omo-epi-cor} it follows that the equation
$f(x)= y$ has at least one solution $x\in U$.
Since the assumption
\[
\sup \{\|k(x)\|: x\in \overline U \}
< \inf \{\|f(x)\|: x\in \partial U \}
\]
implies that $\cc(k(\overline U)) \subseteq B_r$, we get
$X_1 \neq \emptyset$.
Arguing by induction, it is not difficult to show that
$X_n \neq \emptyset$ for any $n \geq 0$ and that the sequence
$\{X_n\}$ is decreasing.
Furthermore, we get
\[
\alpha(X_{n+1})=\alpha(f\sp{-1}(\cc(k(X_n))))
\leq \frac{1}{\omega(f)}\alpha(k(X_n))
\leq \frac{\alpha(k)}{\omega(f)}\alpha(X_n)
\quad \mbox{ for any $n\geq 0$.}
\]
The assumption $\alpha(k)/\omega(f)<1$
implies $\alpha(X_n)\to 0$ as $n\to\infty$.
By the Kuratowski Intersection Theorem (see~\cite{Ku1}), it follows
that the set
\[
X_\infty=\bigcap_{n\geq 0} X_n
\]
is nonempty and compact.
Moreover, it has the following property:
\[
X_\infty \supseteq f\sp{-1}(\cc(k(X_\infty))).
\]
Indeed, notice that for any $n\geq 0$ we have
\[
f\sp{-1}(\cc(k(X_\infty))) \subseteq f\sp{-1}(\cc(k(X_n)))
= X_{n+1}.
\]

By the Dugundji extension theorem (see~\cite{Du}), the restriction
\[
k|_{X_\infty}:X_\infty \to \cc(k(X_\infty))
\]
has a continuous extension
\[
\widetilde k : \overline U  \to \cc(k(X_\infty)),
\]
which is clearly a compact map.
{}From the assumption
\[
\sup \{\|k(x)\|: x\in \overline U \}
< \inf \{\|f(x)\|: x\in \partial U \}
\]
it follows that $\| \widetilde k(x)\| < \|f(x)\|$
for any $x \in \partial U$.
Hence, Corollary~\ref{omo-epi-cor} implies that
the equation
$f(x)=\widetilde k (x)$ has a solution $\bar x \in U$.
Notice that $f( \bar x) \in \cc (k(X_\infty))$
since $\widetilde k (\overline U) \subseteq  \cc(k(X_\infty))$.
Consequently,
$\bar x \in f\sp{-1}(\cc(k(X_\infty)))\subseteq X_\infty$.
Thus, $\widetilde k (\bar x) = k (\bar x)$, which implies
$f(\bar x)=k(\bar x)$.
This completes the proof.
\end{proof}

\medskip

The following result, which is a consequence of the above coincidence
theorem, shows that when $\omega(f)>0$
the two conditions ``$f$ is zero-epi on $U$'' and
``$f$ is not zero-epi on $U$'' are stable under small perturbations.

\begin{theorem}[Stability theorem for zero-epi maps]
\label{teo-stabilita}
Let $U$ be a bounded open subset of $E$,
$f\in C(\overline U,F)$ a $0$-admissible map with $\omega(f)>0$,
and $k\in C(\overline U,F)$.
The following assertions hold:
\begin{enumerate}
\item if $f$ is not $0$-epi, then so is $f+k$ provided that
$\sup \left\{\|k(x)\|: x\in \overline U \right\}$
is sufficiently small; \smallskip
\item if $f$ is $0$-epi, then so is $f+k$ provided that
$\alpha(k)$ and $\sup \left\{\|k(x)\|: x\in \overline U \right\}$
are sufficiently small.
\end{enumerate}
\end{theorem}

\begin{proof}
(1)\
Assume that $f$ is not $0$-epi.
Then, there exists a compact map
$h:\overline U \to F$ such that $h(x)=0$ for all $x \in \partial U$
and the equation $f(x)=h(x)$ has no solutions in $U$.
The fact that $\omega(f-h)=\omega(f)>0$ implies that $f-h$ is
proper on $\overline U$.
Consequently, $(f-h)(\overline U)$ is a
closed subset of $F$.
It follows that $\dist(0,(f-h)(\overline U))>0$.
Assume now that
\[
\sup \left\{\|k(x)\|: x\in \overline U \right\} <
\dist(0,(f-h)(\overline U)).
\]
Then, the equation $f(x)+k(x)=h(x)$ has no solutions in $U$.
Hence, $f+k$ is not $0$-epi.

(2)\
Assume that $f$ is $0$-epi
and let $h:\overline U \to F$ be any compact map such that $h(x)=0$
for all $x \in \partial U$.
Clearly, the map $f-h$ is $0$-epi with
$\omega(f-h)=\omega(f)>0$ and
\[
\inf \{\|f(x)-h(x)\|: x\in \partial U \}=
\inf \{\|f(x)\|: x\in \partial U \}>0
\]
since $f$ is $0$-admissible and proper on $\overline U$.
Assume now that $\alpha(k)<\omega(f)$ and
\[
\sup \{\|k(x)\|: x\in \overline U \}
< \inf \{\|f(x)\|: x\in \partial U \}.
\]
Then, by Theorem~\ref{teo-perturba}, the equation
$f(x)+k(x)-h(x)=0$ has a solution in $U$.
Consequently, $f+k$ is $0$-epi.
\end{proof}


\subsection{Zero-epi maps at a point}
We introduce now the following definitions.
Let $U$ be open in $E$, $f\in C(U,F)$ and $p\in U$.

\begin{definition}
Given $y$ in $F$, we say that $f$ is \textit{$y$-admissible at $p$} if
$f(p)=y$ and $f(x)\neq y$ for any $x$ in a pinched neighborhood of $p$.
\end{definition}

Notice that the map $f$ is $y$-admissible at $p$ if and only if
$f(p)=y$ and $f_p$ is $0$-admissible at $0$.
Furthermore, observe that if $f$ verifies $d_p(f)>0$ then it is
$f(p)$-admissible at $p$.

\begin{definition}
We say that $f$ is \textit{$y$-epi at $p$} if it is $y$-admissible at
$p$ and $y$-epi on any sufficiently small neighborhood of $p$.
\end{definition}

\begin{remark}
In view of the localization property of zero-epi maps
(see Proposition~\ref{prop-local1}),
in the previous definition it is not restrictive to require that there
exists a bounded open neighborhood $U$ of $p$ 
such that $f(x) \neq y$
for all $x \in \overline U$, $x \neq p$, and $f$ is $y$-epi on $U$.
\end{remark}

Notice that $f$ is $y$-epi at $p$ if and only if $f(p)=y$ and $f_p$ is
$0$-epi at $0$.

\medskip

The following local surjectivity property can be deduced from the
corresponding property of zero-epi maps
(see Theorem~\ref{epi-onto} above).

\begin{corollary} \label{cor-onto}
Let $f\in C(U,F)$ be $y$-epi at $p$ and locally proper at $p$.
Then, $y$ is an interior point of the image $f(U)$.
\end{corollary}

\medskip

The following properties are straightforward consequences of the definition.

\begin{proposition} \label{prop-omeoloc}
If $f\in C(U,F)$ is a local homeomorphism at $p$,
then it is $f(p)$-epi at $p$.
\end{proposition}

\begin{proposition} \label{cor-cf}
Let $f\in C(U,F)$ be $y$-epi at $p$, and $c \in \K \backslash \{0\}$.
Then, the map $cf$ is $cy$-epi at $p$.
\end{proposition}

\begin{proposition} \label{cor-Af}
Let $f\in C(U,F)$ be $y$-epi at $p$, and $A:E \to E$,
$B:F \to F$ be linear isomorphisms.
Then, $fA$ is $y$-epi at $A\sp{-1}p$, and $Bf$ is $By$-epi at $p$.
\end{proposition}

\medskip

We observe that a bounded linear operator $L: E\to F$ is $y$-admissible
at $p$ if and only if $Lp=y$ and $L$ is injective.
Moreover, if $Lp=y$, then $L$ is $y$-epi at $p$ if and only if it is
$0$-epi at $0$.

Let $L: E \to F$ be a linear isomorphism.
As a consequence of the Schauder Fixed Point Theorem, it is not
difficult to prove that $L$ is zero-epi on any bounded neighborhood of
the origin (see~\cite{fumavi80}).
In particular, this implies that $L$ is $0$-epi at $0$.


\section{Regular maps and the spectrum at a point}

In this section we define the spectrum of a map at a point.
We need first to introduce the notion of regular map at a
point.


\subsection{Regular maps at a point}

Let $U$ be an open subset of $E$, $f\in C(U,F)$,
and $p\in U$.

\begin{definition}
The map $f$ is said to be \textit{regular at $p$} if the following
conditions hold:
\begin{itemize}
\item[i)]   $d_p(f)>0$; \smallskip
\item[ii)]  $\omega_p(f)>0$; \smallskip
\item[iii)] $f_p$ is $0$-epi at $0$.
\end{itemize}
\end{definition}

Notice that $f$ is regular at $p$ if and only if and $f_p$ is regular
at $0$.
Moreover, if $f$ is regular at $p$ and $c \neq 0$, then $cf$ is
regular at $p$ as well.

\medskip

The following stability property for regular maps,
which can be regarded as a Rouch\'e--type theorem,
will be used in the sequel.

\begin{theorem}[Stability theorem for regular maps]
\label{teo-utile}
Assume that $f$ is regular at $p$ and let $g=f+k$, where $k\in C(U,F)$
is such that $\alpha_p(k)<\omega_p(f)$ and $|k|_p<d_p(f)$.
Then $g$ is regular at $p$.
\end{theorem}

\begin{proof}
By the properties of $d_p$ and $\omega_p$, we have $d_p(g)>0$
and $\omega_p(g)>0$. In particular, the map $g_p$ is $0$-admissible
at $0$. Hence, it is enough to show that $g_p$ is $0$-epi at $0$.

To see this, without loss of generality, we may assume $p=0\in U$
and $f(0)=g(0)=0$. It is possible to choose an open ball
$B_r\subseteq U$ such that the following properties hold:
\begin{itemize}
\item[i)] $\omega(f|_{\overline B_r})>\alpha(k|_{\overline B_r})$;
\smallskip
\item[ii)] there exist positive constants $c$ and $\eta$, with
$c-\eta>0$, such that
\[
\frac{\|f(x)\|}{\|x\|}\geq c \quad \mbox{and} \quad
\frac{\|k(x)\|}{\|x\|}\leq \eta, \quad \mbox{ for any } x\in
\overline B_r;
\]
\item[iii)] $f$ is $0$-epi on $B_r$.
\end{itemize}
Let now $h\in C(\overline B_r,F)$ be compact with $h(x)=0$ for any
$x\in \partial B_r$.
Let us show that the equation $g(x)=h(x)$ has a solution in $B_r$.
Notice that the map $f-h$ is $0$-epi on $B_r$ and
\[
\omega((f-h)|_{\overline B_r})=
\omega(f|_{\overline B_r})>\alpha(k|_{\overline B_r}).
\]
Moreover,
\[
\sup \{\|k(x)\|: x \in \overline B_r \} \leq \eta r < cr \leq
\inf \{\|f(x)\|: x \in \partial B_r \} =
\inf \{\|f(x)-h(x)\|: x \in \partial B_r \}.
\]
Hence, Theorem~\ref{teo-perturba} implies that the equation
$f(x)+k(x)=h(x)$ has a solution in $B_r$.
That is, $g$ is $0$-epi on $B_r$.
Thus, $g$ is $0$-epi at $0$, completing the proof.
\end{proof}

\medskip

Notice that a bounded linear operator $L: E\to F$ is regular at
$p$ if and only if it is regular at $0$. The following proposition
characterizes the bounded linear operators which are regular at
$0$.

\begin{proposition} \label{prop-reg-lineari}
Let $L: E\to F$ be a bounded linear operator.
Then $L$ is regular at $0$ if and only if it is an isomorphism.
\end{proposition}

\begin{proof}
Let $L$ be an isomorphism.
Then, as already pointed out, the Schauder Fixed Point Theorem
implies that $L$ is $0$-epi at $0$.
Moreover, $d(L)>0$ by
Proposition~\ref{prop-lin1-dq}-(3), and
$\omega(L)\geq d(L)>0$ by
Proposition~\ref{prop-lin1-dq}-(2).

Conversely, assume that $L$ is regular at $0$.
Then $L$ is injective since $d(L)>0$, and locally proper at $0$ by
Proposition~\ref{prop-alfaomega}-(6), as $\omega(L)>0$.
Hence, $L$ being $0$-epi at $0$, Corollary~\ref{cor-onto} implies that
$L$ is surjective as well.
\end{proof}


\subsection{Spectrum at a point}

Let now $f\in C(U,E)$ and $p\in U$.
We define the
\textit{spectrum of the map $f$ at the point $p$} as the set
\[
\sigma(f,p)=\{\lambda \in \K: \lambda-f \mbox{ is not regular at } p \},
\]
where $\lambda - f$ stands for $\lambda I - f$, $I$ being the identity
on $E$.

It is convenient to define the following subset of $\sigma(f,p)$:
\[
\sigma_\pi(f,p)=\{\lambda \in \K: d_p(\lambda - f)=0 \mbox{ or }
\omega_p(\lambda - f)=0\}.
\]
We shall call
$\sigma_\pi (f,p)$ the
\textit{approximate point spectrum of $f$ at $p$}.

Clearly, for a bounded linear operator $L: E\to E$,
$\sigma(L,p)$ and $\sigma_\pi(L,p)$ do not depend on $p$.
Hence, we can simply write $\sigma(L)$ and $\sigma_\pi(L)$ instead of
$\sigma(L,p)$ and $\sigma_\pi(L,p)$.

The above notation and definitions are justified by the following
result, which is the analogue of Theorem~8.1.1 in~\cite{fumavi78}.

\begin{theorem} \label{teo-spettro-lineari}
Let $L: E\to E$ be a bounded linear operator.
Then
\begin{enumerate}
\item $\sigma(L)$ coincides with the usual spectrum of $L$; \smallskip
\item $\sigma_\pi(L)$ coincides with the usual approximate point
spectrum of $L$.
In other words, $\lambda \in \sigma_\pi (L)$ if and only if
\begin{equation} \label{eq-aps}
\inf_{\|x\|=1} \|\lambda x - Lx\|=0.
\end{equation}
\end{enumerate}
\end{theorem}

\begin{proof}
(1)\
The assertion follows from Proposition~\ref{prop-reg-lineari}.

(2)\
Let $\lambda \in \K$ be such that \eqref{eq-aps} holds.
By Remark~\ref{oss-dq-posom}, we have $d(\lambda -L)=0$.
Thus, $\lambda \in \sigma_\pi (L)$.
Conversely, assume that $\lambda \in \sigma_\pi (L)$.
Since Proposition~\ref{prop-lin1-dq}-(2) implies
$\omega(\lambda -L) \geq d(\lambda -L)$, we may assume
$d(\lambda -L)=0$.
Hence,
\[
\inf_{\|x\|=1} \|\lambda x - Lx\|=0
\]
and the assertion follows.
\end{proof}

\medskip

The next result, which is a consequence of
Theorem~\ref{teo-utile}, provides a sufficient condition for the
spectrum of $f$ at $p$ to be bounded.
Set
\[
q_p(f)=\max\{\alpha_p(f),|f|_p\}
\]
and define the \textit{spectral radius of $f$ at $p$} as
\[
r_p(f)=\sup\{|\lambda|:\lambda\in\sigma(f,p)\}.
\]

\begin{proposition} \label{prop-limitato}
We have $r_p(f) \leq q_p (f)$.
In particular, if $f$ is locally compact and
quasibounded at $p$, then $\sigma (f,p)$ is bounded.
\end{proposition}

\begin{proof}
We may assume that $q_p(f)< + \infty$.
Notice that, if $\lambda \neq 0$, then $\lambda I$ is regular at $p$.
By the properties of $d_p$ and $\omega_p$ we get
$\omega_p(\lambda I) = |\lambda| = d_p(\lambda I)$.
Hence, Theorem~\ref{teo-utile} implies that $\lambda - f$ is regular
at $p$ when $|\lambda| > q_p$.
\end{proof}

In the case of a bounded linear operator $L : E \to E$, the number
$r_p(L)$ is independent of $p$.
Therefore, this number will be denoted by $r(L)$.

{}From Proposition~\ref{prop-limitato} above
(and Proposition~\ref{prop-lin1-dq}-(1))
we recover the well known property
\[
r(L) \leq \|L\|.
\]

\medskip

The following, which is the analogue of Theorem~8.1.2
in~\cite{fumavi78}, is our main result.

\begin{theorem} \label{teoremone-sigma}
The following properties hold.

\begin{enumerate}
\item $\sigma(f,p)$ is closed; \smallskip
\item $\sigma_\pi(f,p)$ is closed; \smallskip
\item $\sigma(f,p)\backslash\sigma_\pi(f,p)$ is open; \smallskip
\item $\partial \sigma(f,p)\subseteq\sigma_\pi(f,p)$.
\end{enumerate}
\end{theorem}

\begin{proof}
(1)\
Let $\lambda \not \in \sigma (f,p)$ and $\mu \in \K$ be such that
\[
|\mu| < \min \{ d_p (\lambda-f),\omega_p (\lambda-f) \}.
\]
By Theorem~\ref{teo-utile}, the map $(\lambda+\mu)-f$ is regular at
$p$. Hence, $\K \backslash \sigma(f,p)$ is open.

(2)\
Is a straightforward consequence of the fact that
$d_p (\lambda - f)$ and $\omega_p (\lambda - f)$ depend continuously
on $\lambda$ (see Proposition~\ref{prop-alfaomega} and
Proposition~\ref{prop-dq-proprieta}).

(3)\
Let $\lambda \in \sigma (f,p) \backslash \sigma_\pi (f,p)$.
That is, $d_p (\lambda-f) >0$, $\omega_p (\lambda-f) >0$,
and the map $\lambda-f_p$ is not $0$-epi at $0$.
It is sufficient to show that, if $\mu$ is sufficiently small, the map
$(\lambda+\mu) - f_p$ is not $0$-epi at $0$.

To see this, choose a suitably small open neighborhood $V$ of $0$
such that $\omega ((\lambda-f_p)|_{\overline V}) >0$.
Since the map $\lambda-f_p$ is not $0$-epi on $V$, there exists a
compact map $h:\overline V \to E$ such that $h(x)=0$ for all $x
\in \partial V$ and the equation $\lambda x-f_p(x)=h(x)$ has no
solutions in $V$.
{}From
\[
\omega((\lambda-f_p-h)|_{\overline V})=
\omega((\lambda-f_p)|_{\overline V})>0
\]
it follows that the map $\lambda-f_p-h$ is proper on $\overline V$,
which implies that $(\lambda-f_p-h)(\overline V)$ is closed.
Thus, there exists $\eta>0$ such that
$\|\lambda x-f_p(x)-h(x)\|\geq\eta$ for all $x\in\overline V$.

Now, let $\mu$ be such that
$|\mu| \sup \{ \|x\| : x \in \overline V \} < \eta$.
We have
\[
\|(\lambda+\mu) x-f_p(x)-h(x)\| \geq \eta - |\mu| \|x\|>0
\]
for all $x \in \overline V$.
Hence, $(\lambda+\mu) -f_p-h$ is not $0$-epi on $V$ and
so is $(\lambda+\mu) -f_p$.
Since $V$ is arbitrarily small, the map $(\lambda+\mu) -f_p$ is not
$0$-epi at $0$ on the basis of Proposition~\ref{prop-local1}, and the
assertion follows.

(4)\
Let $\lambda \in \partial \sigma(f,p)$. Since
$\sigma(f,p)$ is closed we get
$\lambda \in \sigma(f,p)$.

By contradiction, assume that $\lambda\not\in\sigma_\pi(f,p)$.
Since the set $\sigma(f,p)\backslash\sigma_\pi(f,p)$ is open in $\K$,
we have $\lambda \in \Int(\sigma(f,p))$, contradicting the assumption
$\lambda\in\partial\sigma(f,p)$.
\end{proof}

\begin{corollary} \label{cor-componenti-sigma}
Let $W$ be a connected component of
$\K \backslash \sigma_\pi (f,p)$.
Then, $W$ is open in $\K$ and maps of the form $\lambda-f_p$, with
$\lambda \in W$, are either all $0$-epi at $0$ or all not $0$-epi at
$0$.
\end{corollary}

\begin{proof}
The set $W$ is open since $\K \backslash \sigma_\pi(f,p)$ is locally
connected.
Let
\[
A=\{\lambda \in W: \lambda - f_p \mbox{ is
   $0$-epi at } 0 \}
\quad \mbox{ and } \quad
B=\{\lambda \in W: \lambda - f_p \mbox{ is not
   $0$-epi at } 0 \}.
\]
Since $W\cap \sigma_\pi(f,p)=\emptyset$,
we have
\[
A=W\backslash\sigma(f,p),
\]
which is open by Theorem~\ref{teoremone-sigma}-(1),
and
\[
B=W\cap(\sigma(f,p)\backslash\sigma_\pi(f,p)),
\]
which is open by Theorem~\ref{teoremone-sigma}-(3).
The connectedness of $W$ implies either
$W=A$ or $W=B$.
\end{proof}

\begin{corollary} \label{cor-sigma-limitato}
The following assertions hold.
\begin{itemize}
\item[i)] Assume that $\K=\C$, $\sigma(f,p)$ is bounded, and $\lambda$
belongs to the unbounded component of $\C \backslash \sigma_\pi(f,p)$.
Then $\lambda - f$ is regular at $p$.
In particular (as pointed out by J.R.L.\ Webb in a private
communication), if $\sigma_\pi(f,p)$ is countable, then
$\sigma(f,p)= \sigma_\pi(f,p)$.
\smallskip
\item[ii)] Assume that $\K=\R$, $\sigma(f,p)$ is bounded from above
(resp.\ below), and $\lambda$ belongs to the right (resp.\ left)
unbounded component of $\R\backslash\sigma_\pi(f,p)$.
Then $\lambda - f$ is regular at $p$.
\end{itemize}
\end{corollary}


\subsection{A finer decomposition of the spectrum at a point}

Given $f$ and $p$ as above, we introduce the following subsets of
$\sigma_\pi (f,p)$:
\[
\sigma_\omega (f,p) =
\{\lambda \in \K: \omega_p (\lambda - f) =0 \}
\quad \mbox{ and } \quad
\Sigma (f,p) =
\{\lambda \in \K: \d_p (\lambda - f)=0 \}.
\]
Evidently, $\sigma_\pi (f,p) = \sigma_\omega (f,p) \cup \Sigma (f,p)$.
We point out that $\Sigma (f,p)$ has been introduced in~\cite{fuvi77}
by the last two authors.

\begin{remark} \label{oss-sigmapi}
If $E$ is finite dimensional, then $\sigma_\omega(f,p)=\emptyset$ for
any $p$ and hence $\sigma_\pi (f,p) = \Sigma (f,p)$.
\end{remark}

\begin{remark} \label{oss-real}
It is interesting to observe that, for a real function $f$, the
spectrum $\sigma(f,p)$ is completely 
determined by the \textit{Dini's derivatives} of $f$ at $p$,
that is, by the following four extended real numbers:
\[
D_{-}f(p) = \liminf_{h \to 0^-} \frac{f(p+h) - f(p)}{h}, \quad
D\sp{-}f(p) = \limsup_{h \to 0^-} \frac{f(p+h) - f(p)}{h},
\]
\[
D_{+}f(p) = \liminf_{h \to 0^+} \frac{f(p+h) - f(p)}{h}, \quad
D\sp{+}f(p) = \limsup_{h \to 0^+} \frac{f(p+h) - f(p)}{h}.
\]
It is not difficult to show that $\sigma(f,p)$ is the closed subinterval
of $\R$ whose endpoints are, respectively, the smallest and
the largest of the Dini's derivatives.
Thus, any closed interval (the empty set, a singleton and $\R$
included) is the spectrum at a point of some continuous function.
For example, if all the four Dini's derivatives of $f$ at
$p$ are $+\infty$ (or $-\infty$), then $\sigma(f,p) = \emptyset$.

We point out that, 
even in the one dimensional real case, $\Sigma(f,p)$ need not coincide
with $\sigma(f,p)$.
In fact, one can check that $\Sigma(f,p)$ is the union of two closed
intervals: one with endpoints 
$D_{-}f(p)$ and $D\sp{-}f(p)$, and the other one
with endpoints $D_{+}f(p)$ and $D\sp{+}f(p)$.
Hence, $\sigma(f,p)$ is the smallest interval containig $\Sigma(f,p)$,
and this agrees with Theorem  \ref{teoremone-sigma}.

As a consequence of these facts,
if $f$ is differentiable at $p$ one gets
$\sigma(f,p) = \Sigma(f,p) = \{ f'(p) \}$.
This agrees with Corollary~\ref{cor-equivalenza} below.
\end{remark}

\medskip

The following simple examples illustrate two ``pathological cases''.

\begin{example}
Let $f:\R \to \R$ be defined by
\[
f(x) = \sign(x) \sqrt{|x|}.
\]
Then, $\Sigma(f,0) = \sigma(f,0) = \emptyset$.
\end{example}

\begin{example}
Let $f:\R \to \R$ be defined by
\[
f(x) = \sqrt{|x|}.
\]
Then, $\Sigma(f,0) = \emptyset$ and $\sigma(f,0) = \R$.
\end{example}

\medskip

The next propositions show some properties of
$\sigma_\omega (f,p)$ and $\Sigma (f,p)$.

\begin{proposition} \label{teo-sigmapi2}
The sets $\sigma_\omega(f,p)$ and $\Sigma(f,p)$ are closed.
\end{proposition}

\begin{proof}
The assertion follows directly from the fact that
$d_p(\lambda - f)$ and $\omega_p(\lambda - f)$ depend
continuously on $\lambda$.
\end{proof}

\begin{proposition} \label{prop-stime}
The following estimates hold:
\begin{enumerate}
\item $\lambda \in \Sigma(f,p)$
implies $d_p(f) \leq |\lambda| \leq |f|_p$; \smallskip
\item $\lambda \in \sigma_\omega (f,p)$
implies $\omega_p(f) \leq |\lambda| \leq \alpha_p(f)$.
\end{enumerate}
\end{proposition}

\begin{proof}
(1)\ Let $\lambda \in \Sigma(f,p)$.
We have $d_p(\lambda I) = |\lambda| = |\lambda I|_p$
by Proposition~\ref{prop-dq-proprieta}-(1).
Therefore, Proposition~\ref{prop-dq-proprieta}-(5) implies
$0 = d_p(\lambda-f) \geq d_p(f)-|\lambda|$ and
$0 = d_p(\lambda-f) \geq |\lambda|-|f|_p$.
Consequently, $d_p(f) \leq |\lambda| \leq |f|_p$.

(2)\ Let $\lambda \in \sigma_\omega (f,p)$.
We have $\omega_p(\lambda I) = |\lambda| = \alpha_p(\lambda I)$
by Proposition~\ref{prop-alfaomega}-(1).
Thus, by Proposition~\ref{prop-alfaomega}-(4), we get
$0 = \omega_p(\lambda-f) \geq \omega_p(f)-|\lambda|$ and
$0 = \omega_p(\lambda-f) \geq |\lambda|- \alpha_p(f)$.
Hence, $\omega_p(f) \leq |\lambda| \leq \alpha_p(f)$.
\end{proof}

\begin{proposition} \label{prop-spettro-compatte}
Assume $\dim E = + \infty$ and let $f$ be locally compact.
Then, $0 \in \sigma(f,p)$.
More precisely, $\sigma_\omega(f,p) = \{ 0 \}$.
Thus, $\sigma_\pi(f,p)= \{0\} \cup \Sigma(f,p)$.
\end{proposition}

\begin{proof}
Since $f$ is locally compact, Proposition~\ref{prop-alfaomega}-(5)
implies $\alpha_p (f)=0$.
Consequently, by Proposition~\ref{prop-alfaomega}, we get
$\omega_p (\lambda-f) = \omega_p (\lambda I) = |\lambda|$
and the assertion follows.
\end{proof}

{}From Proposition~\ref{prop-stime} it follows that, when $|f|_p$ and
$\alpha_p (f)$ are finite, then $\sigma_\pi (f,p)$ is bounded.
Actually, Proposition~\ref{prop-limitato} shows that, in this case,
the whole spectrum $\sigma (f,p)$ is bounded.

\medskip

The following is a pathological example of a continuous map
$f:\R \to \R$ with an unbounded approximate point spectrum at a point.

\begin{example}
Let $f:\R \to \R$ be defined by
\[
f(x) =
\left\{
\begin{array}{ll}
\displaystyle \sqrt{|x|} \sin \frac 1x  & x \neq 0,\\
0 & x=0.
\end{array}
\right.
\]
Clearly, $\sigma_\pi(f,0) = \Sigma(f,0) = \R$.
\end{example}

\medskip

For a bounded linear operator $L: E\to E$, the sets
$\sigma_\omega(L,p)$ and $\Sigma(L,p)$ do not depend on the point $p$
and will be denoted by $\sigma_\omega(L)$ and $\Sigma(L)$,
respectively.
We remark that this notation agrees with the one introduced
in~\cite{fumavi78}.
In fact, we recall the following definitions.

As in~\cite{fumavi78}, given $f:E \to E$, define
the \textit{asymptotic approximate point spectrum} of $f$ by
\[
\sigma_\pi (f) =\{ \lambda \in \K: d (\lambda - f) =0
\mbox{ or } \omega (\lambda - f) =0 \}
\]
and define
\[
\sigma_\omega (f) =\{ \lambda \in \K: \omega (\lambda - f) =0 \}
\quad \mbox{ and } \quad
\Sigma (f) =\{ \lambda \in \K: d (\lambda - f) =0 \}.
\]

The next result shows that, for a positively homogeneous map
$f:E \to E$, the approximate point spectrum of $f$ at $0$
coincides with the asymptotic approximate point spectrum of $f$.
In particular, the same is true for bounded linear operators.

\begin{proposition} \label{prop-spettro-posom}
Let $f:E \to E$ be positively homogeneous.
Then, $\sigma_\pi (f,0) = \sigma_\pi (f)$.
More precisely,
$\sigma_\omega (f,0) = \sigma_\omega (f)$ and
$\Sigma (f,0)= \Sigma (f)$.
In addition, $\lambda \in \Sigma (f)$ if and only if
\[
\inf_{\|x\|=1} \|\lambda x - f(x)\| =0.
\]
\end{proposition}

\begin{proof}
Is a direct consequence of
Proposition~\ref{prop-alfaposom} and
Remark~\ref{oss-dq-posom}.
\end{proof}

\begin{corollary}
Let $L:E \to E$ be a bounded linear operator.
Then, $\sigma_\pi (L) = \Sigma (L)$.
\end{corollary}

\begin{proof}
The assertion follows from Proposition~\ref{prop-spettro-posom} above and
Theorem~\ref{teo-spettro-lineari}.
\end{proof}

One could show that, for a positively homogeneous map $f$,
the spectrum of $f$ at $0$ coincides with the
\textit{asymptotic spectrum} of $f$ defined in~\cite{fumavi78}.

\medskip

In the case of positively homogeneous maps, it is meaningful to
introduce the concept of eigenvalue.
Let $f:E \to E$ be positively homogeneous.
As in the linear case, we say that $\lambda \in \K$ is an
\textit{eigenvalue of $f$ at $0$} if the equation
$\lambda x = f(x)$
admits a nontrivial solution.

\begin{proposition} \label{prop-autovalori}
Assume $f:E \to E$ positively homogeneous
and $\lambda \not \in \sigma_\omega (f,0)$.
Then, $\lambda \in \Sigma(f,0)$
if and only if $\lambda$ is an eigenvalue of $f$ at $0$.
\end{proposition}

\begin{proof}
Let $\lambda \not \in \sigma_\omega (f,0)$;
that is, $\omega (\lambda - f) >0$.
Then, Proposition~\ref{prop-betaposom} implies that
$\lambda \in \Sigma (f,0)$ if and only if
there exists $x \in E$, $x \neq 0$, such that
$\lambda x = f(x)$; that is, if and only if $\lambda$ is an eigenvalue
of $f$ at $0$.
\end{proof}

\begin{corollary} \label{cor-autovalori}
The following assertions hold.
\begin{itemize}
\item[i)] Assume $\dim E= +\infty$ and let $f:E \to E$ be positively
homogeneous and locally compact.
Then,
$\Sigma(f,0) \backslash \{0\} =
\{ \lambda \in \K \backslash \{0\} :
\lambda \mbox{ eigenvalue of } f \mbox{ at } 0 \}$.
\item[ii)] Assume $E$ finite dimensional and $f:E \to E$ positively
homogeneous. Then,
$\Sigma(f,0) = \{ \lambda \in \K :
\lambda \mbox{ eigenvalue of } f \mbox{ at } 0 \}$.
\end{itemize}
\end{corollary}

\medskip

The next result concerns compact linear operators
(see also~\cite{fuvi75}).

\begin{proposition} \label{prop-spettro-compatte+lineari}
Assume $\dim E= +\infty$ and let $L: E \to E$ be a compact linear
operator. Then, $0\in\Sigma(L)$ and $\Sigma(L) = \sigma(L)$.
\end{proposition}

\begin{proof}
Assume that $0\not\in\Sigma(L)$.
Then, there is a positive constant $m$ such that
\[
\inf_{\|x\|=1} \|L(x)\| \geq m>0.
\]
Hence, $L$ is a linear isomorphism between $E$ and the closed subspace
$E_0 = \Im L$ of $E$.
This is impossible since $\dim E_0 = \dim E = +\infty$ and $L$ is a
compact operator. Consequently, $0\in\Sigma(L)$.

Since $\Sigma(L) \subseteq \sigma(L)$, to prove that
$\Sigma(L) = \sigma(L)$ it is enough to show the opposite inclusion.
Let $\lambda \neq 0$ be given.
Recall that, since $L$ is a compact linear operator,
from the Fredholm alternative it follows that if
$\lambda - L$ is injective then $\lambda - L$ is invertible.
That is, when $\lambda \neq 0$,
the condition $\lambda \not \in \Sigma(L)$ implies
$\lambda \not \in \sigma (L)$.
Now the equality $\sigma(L)=\Sigma(L)$
follows taking into account that $0\in\Sigma(L)$.
\end{proof}

\medskip

Proposition~\ref{prop-componente-sigma} below extends Theorem~2.1
in~\cite{fuvi77} by replacing the assumption
``$f$ quasibounded at $p$'' with the weaker condition
``$\sigma (f,p)$ bounded''.
We wish to stress the fact that, contrary to the proof
of~\cite[Theorem~2.1]{fuvi77}, here no degree theory is involved.
We wish further to observe that a result of this type could not be
stated in~\cite{fuvi77} because of the lack of an exhaustive
notion of spectrum at a point.
The same proposition exhibits an exclusively nonlinear phenomenon since,
in the compact linear case, zero always belongs to the approximate
point spectrum. An example illustrating this peculiarity will be given in
Section \ref{sect-examples}.

\begin{proposition} \label{prop-componente-sigma}
Let $U$ be an open subset of $E$, with $\dim E= +\infty$.
Let $f \in C(U,E)$ be locally compact and assume that
$\sigma (f,p)$ is bounded.
Then, $0 \not \in \Sigma (f,p)$ implies that $0$ is an interior point
of $\sigma (f,p)$.
In particular, the connected component
of $\K \backslash \Sigma (f,p)$ containing $0$ is bounded.
\end{proposition}

\begin{proof}
Without loss of generality, we may assume $p=0 \in U$ and
$f(0)=0$.

Note that the assumption $0 \not \in \Sigma (f,0)$ is equivalent to
$d_0(f) >0$, and this implies that there exist $r>0$, $2c>0$ such that
\[
\|f(x)\| > 2c \|x\|
\]
for all $x \in B_r$. Hence,
\[
\|\lambda x - f(x)\| \geq \|f(x)\| - |\lambda| \|x\| >
(2c-|\lambda|) \|x\| \geq c \|x\|
\]
for all $x \in B_r$ and $\lambda \in \K$ such that
$|\lambda| \leq c$.

We claim that, if $|\lambda| >0$ is sufficiently small, then
the map $\lambda -f$ is not $0$-epi at $0$.

By contradiction, assume that
there exists a sequence $\{ \lambda_n \}$ in
$\K \backslash \{0\}$ such that $\lambda_n \to 0$ as $n \to \infty$
and $\lambda_n -f$ is $0$-epi at $0$ for any $n \in \N$.
Without loss of generality, we may assume
$0< |\lambda_n| \leq c$ for any $n \in \N$.
Furthermore, in view of the antilocalization property of zero-epi maps
(see Proposition~\ref{prop-taglio}),
it is not restrictive to assume that $\lambda_n -f$ is
$0$-epi on $B_r$ for any $n$.
Consequently,
since $\|\lambda_n x - f(x)\| \geq rc$ for all $x \in \partial B_r$
and $n \in \N$, Corollary~\ref{omo-epi-cor} implies that
the equation $\lambda_n x - f(x)=y$ has a solution in $B_r$ for all
$y \in E$ with  $\|y\| < cr$ and all $n \in \N$.

Fix now $y\in E$ with $\|y\| < cr$.
Then, there exists a sequence
$\{x_n\}$ in $B_r$ such that $\lambda_n x_n - f(x_n) =y$.
Since $\{x_n\}$ is bounded and $\lambda_n \to 0$, it follows
that $f(x_n) \to -y$.
Thus, the image $f(B_r)$ is dense in $B_{cr}$.
This is a contradiction since $r>0$ can taken so small that
$\overline{f(B_r)}$ is a compact set which, in this case, has empty
interior (recall that $\dim E = + \infty$).
Hence, $\lambda -f$ is not $0$-epi at $0$ for $|\lambda| >0$
sufficiently small, as claimed.
Therefore, $0$ is an interior point of $\sigma (f,0)$, since
$0 \in \sigma_\omega (f,0)$ on the basis of
Proposition~\ref{prop-spettro-compatte}.

Finally,
the last assertion follows from the assumption that $\sigma (f,0)$ is
bounded, and from the fact that
\[
\partial \sigma (f,p) \subseteq \sigma_\pi (f,0) =
\{0\} \cup \Sigma(f,0)
\]
since $f$ is locally compact
(see Proposition~\ref{prop-spettro-compatte}
and Corollary~\ref{cor-componenti-sigma}).
This completes the proof.
\end{proof}

\medskip

In what follows the notation $\sigma(f,p)\equiv\sigma(g,p)$ stands
for $\sigma(f,p)=\sigma(g,p)$,
$\sigma_\omega(f,p)=\sigma_\omega(g,p)$, and
$\Sigma(f,p)=\Sigma(g,p)$.
Recall that $q_p(f)=\max\{\alpha_p(f),|f|_p\}$.

\begin{theorem} \label{teo-decomp-spettrale}
Given an open subset $U$ of $E$, $f,g \in C(U,E)$ and
$p\in U$, one has
\begin{enumerate}
\item $\sigma(c f,p) \equiv c \sigma (f,p)$, for any
   $c \in \K$;\smallskip
\item $\sigma(c + f,p) \equiv c + \sigma (f,p)$, for any
   $c \in \K$;\smallskip
\item $q_p (f - g)=0$ implies $\sigma(f,p)\equiv\sigma(g,p)$.
\end{enumerate}
\end{theorem}

\begin{proof}
(1)\
The equality is trivial if $c=0$. Assume $c \neq 0$.
One has $\sigma_\omega(c f,p)=c \sigma_\omega(f,p)$ and
$\Sigma(c f,p)=c \Sigma(f,p)$.
By Proposition~\ref{cor-cf},
given $\lambda \in \K$, $\lambda - f_p$ is $0$-epi at $0$ if
and only if so is $c(\lambda - f_p)$.
Thus, $\sigma(c f,p)=c \sigma(f,p)$ and the assertion follows.

(2)\ Analogous to (1).

(3)\
{}From $q_p (f - g)=0$ it follows that $f=g+h$, where
$\alpha_p (h) = |h|_p=0$. In particular, $|h|_p=0$ implies that
$\Sigma (f,p) = \Sigma (g,p)$, since
$d_p (\lambda-f) = d_p (\lambda-g)$ for any $\lambda \in \K$
by Proposition~\ref{prop-dq-proprieta}-(5).
Moreover, from $\alpha_p (h) =0$
it follows that $\sigma_\omega (f,p) = \sigma_\omega (g,p)$, since
$\omega_p (\lambda-f) = \omega_p (\lambda-g)$ for any
$\lambda \in \K$ by Proposition~\ref{prop-alfaomega}-(4).
Finally, let $\lambda \in \K$ be given.
Since $\alpha_p (h) =|h|_p=0$,
Theorem~\ref{teo-utile} implies that $\lambda-f$ is regular at $p$ if
and only if so is $\lambda-g$.
This shows that $\sigma (f,p) = \sigma (g,p)$, and the proof is
complete.
\end{proof}

The following result extends a well known property of the spectrum of
a linear operator.

\begin{theorem}
Let $U$ be an open subset of $E$, $f \in C(U,E)$ and
$p \in U$.
Given a linear isomorphism $A:E \to E$, let $q=Ap$.
Then,
\[
\sigma(f,p) \equiv \sigma (AfA\sp{-1},q).
\]
\end{theorem}

\begin{proof}
The equality $\Sigma (f,p) = \Sigma (AfA\sp{-1},q)$ has been
established in~\cite[Proposition~2.2]{fuvi77}.

Let us show that
$\sigma_\omega (f,p) = \sigma_\omega (AfA\sp{-1},q)$.
By Remark~\ref{oss-sigmapi}, we may assume $\dim E = +\infty$.
Let $\lambda \in \K$ be given.
Observe that $\lambda -AfA\sp{-1} = A(\lambda-f)A\sp{-1}$.
Thus, it is enough to prove that $\omega_p (\lambda-f) =0$ if and only
if $\omega_q (A(\lambda-f)A\sp{-1}) =0$.
Assume $\omega_q (A(\lambda-f)A\sp{-1}) =0$.
By Proposition~\ref{prop-alfacomposto} we get
\[
\omega_q (A(\lambda-f)A\sp{-1}) \geq
\omega (A) \,\omega_p (\lambda-f) \,\omega (A\sp{-1}),
\]
and Proposition~\ref{prop-alfaomega}-(7) implies
that $\omega (A)$ and $\omega (A\sp{-1})$
are different from zero.
Thus, $\omega_p(f)=0$.
The converse implication can be proved in an analogous way.

Finally, given $\lambda \in \K$,
Proposition~\ref{cor-Af} implies that
$\lambda - f_p$ is $0$-epi at $0$ if and only if so is
$A(\lambda - f)A\sp{-1}_q$.
Hence, $\sigma (f,p) = \sigma (AfA\sp{-1},q)$, and the assertion
follows.
\end{proof}

\medskip

For $C \sp 1$ maps we have the following result which is a
direct consequence of
Theorem~\ref{teo-decomp-spettrale}-(3).

\begin{corollary} \label{cor-equivalenza}
Let $f:U \to E$ be of class $C \sp 1$ and $p \in U$.
Then, $\sigma (f,p) \equiv \sigma (f'(p))$.
\end{corollary}

Notice that the equivalence $\sigma (f,p) \equiv \sigma (f'(p))$ holds
true even when the map $f \in C(U,E)$ is merely Fr\'echet
differentiable at the point $p \in U$, provided that the remainder
$\phi \in C(U,E)$, defined
as
\[
\phi (x)= f(x) - f'(p)(x-p), \quad x \in U,
\]
verifies $\alpha_p (\phi) =0$.
This is the case, for instance, if $f=g+h$, where $g$ is of class
$C \sp1$ and $h$ is locally compact and Fr\'echet differentiable at $p$
(but not necessarily $C \sp1$).
As an example, consider the map $f:E \to E$ defined by
\[
f(x)= \left\{
\begin{array}{ll}
x + \|x\|\sp 2 \left( \sin \frac{1}{\|x\|} \right) v
& \mbox{if } x \neq 0,\\
0 & \mbox{if } x=0,
\end{array}
\right.
\]
where $v \in E \backslash \{0\}$ is given, and $p=0$.

\medskip

In the case when the map $f$ is of class $C \sp 1$,
Corollary~\ref{cor-equivalenza} above implies that the
multivalued map that associates to every point $p$ the spectrum
$\sigma (f,p) \subseteq \K$ is upper semicontinuous.
This depends on the well known fact that
so is the map which associates to any bounded linear operator
$L:E \to E$ its spectrum $\sigma (L)$.
The next one-dimensional example shows that the multivalued map
$p \multimap \sigma (f,p) \subseteq \K$ need not be semicontinuous
if $f$ is merely $C \sp 0$.

\begin{example}
Let $f:\R \to \R$ be defined by
\[
f(x) = x^2 \sin \frac{1}{x}.
\]
Notice that $f$ is $C \sp 1$ on $\R \setminus \{0\}$ and
merely differentiable at $0$ with $f'(0)=0$.
Thus, as pointed out above, $\sigma (f,0) = \{ 0 \}$.
Consequently, the map $p \multimap \sigma (f,p) = \{ f'(p) \}$
is not upper semicontinuous at $0$ 
since $f'$ is not continuous at $0$.
\end{example}

\medskip

The following property of nonemptiness of the spectrum at a point is a
straightforward consequence of Corollary~\ref{cor-equivalenza} above.

\begin{corollary}
If $\K=\C$ and $f:U\to E$ is of class $C\sp 1$, then for any $p\in U$ the
spectrum $\sigma(f,p)$ is nonempty.
\end{corollary}

We give now an example of a complex map with empty spectrum.

\begin{example}
Define $f: \C\sp 2 \to \C\sp 2$ by
\[
f(z,w) = (\bar w, i \bar z), \quad (z,w)\in \C\sp2.
\]
By Remark~\ref{oss-sigmapi}, we have
$\sigma_\pi(f,0) = \Sigma(f,0)$.
Moreover, since $f$ is positively homogeneous,
by Proposition~\ref{prop-spettro-posom}
we have $\Sigma(f,0) = \Sigma(f)$,
and the fact that $\Sigma (f) = \emptyset$ has been
established in~\cite{gema} (see also~\cite{fumavi78}).
Now,
to prove that $\sigma (f,0) = \emptyset$, recall the following two
facts:
\begin{itemize}
\item[(i)]  $\partial \sigma (f,0) \subseteq \sigma_\pi (f,0)$;
\smallskip
\item[(ii)] $\lambda \in \K$, $|\lambda|>q_0(f)$ implies that
$\lambda - f$ is regular at $0$.
\end{itemize}
Since, in our case, $\sigma_\pi (f,0)$ is empty, then either
$\sigma(f,0)=\C$ or $\sigma(f,0)=\emptyset$.
Now, $q_0(f)=|f|_0=1$ implies that $\sigma(f,0)\not= \C$.
Therefore $\sigma(f,0)=\emptyset$.
\end{example}


\section{Bifurcation in the non-differentiable case}

Let $U$ be an open subset of $E$ and $f \in C(U,E)$.
Assume that $0\in U$ and $f(0)=0$, and consider the equation
\begin{equation} \label{equ-bif}
\lambda x =f(x), \quad \lambda \in \K.
\end{equation}
A solution $(\lambda,x)$ of~\eqref{equ-bif} is called
\textit{nontrivial} if $x\neq 0$.
We recall that $\lambda \in \K$ is a \textit{bifurcation point} for
$f$ if any neighborhood of $(\lambda,0)$ in $\K \times E$ contains a
nontrivial solution of~\eqref{equ-bif}.
We will denote by $\B (f,0)$ the set of bifurcation points of $f$.
Notice that $\B (f,0)$ is closed since $\lambda \in \B(f,0)$ if and
only if $(\lambda,0)$ belongs to the closure $\overline S$ of the set
$S$ of the nontrivial solutions of~\eqref{equ-bif}.

It is well known that if $f$ is Fr\'echet differentiable at $0$ then
the set $\B (f,0)$ of bifurcation points of $f$ is contained in the
spectrum $\sigma (f'(0))$ of the Fr\'echet derivative $f'(0)$.
The next proposition (see also~\cite{fuvi77}) extends this necessary
condition. The simple proof is given for the sake of completeness.

\begin{proposition} \label{prop-bif}
The set $\B(f,0)$ is contained in $\Sigma(f,0)$,
and hence in $\sigma_\pi(f,0)$.
\end{proposition}

\begin{proof}
Let $\lambda\in\B(f,0)$.
Then, there exists a sequence $\{(\lambda_n,x_n)\}$, with $x_n\neq 0$
and $\lambda_n x_n=f(x_n)$ for all $n\in\N$,
such that $(\lambda_n,x_n) \to (\lambda,x)$.
Thus,
\[
\frac{\|\lambda x_n - f(x_n)\|}{\|x_n\|} =
\frac{\|\lambda x_n - \lambda_n x_n\|}{\|x_n\|} =
|\lambda_n-\lambda|.
\]
This shows that $\lambda\in\Sigma(f,0)$.
\end{proof}

\begin{remark} \label{oss-autovalori}
As in the linear case, for a positively homogeneous map $f:E \to E$ we
have that if $\lambda \in \K$ is an eigenvalue of $f$ at $0$ then
$\lambda\in\B(f,0)$.
If, moreover, $\omega (\lambda - f) >0$, the converse is also true in
view of Proposition~\ref{prop-autovalori} and
Proposition~\ref{prop-bif} above.
\end{remark}

The following result provides a sufficient condition for the
existence of bifurcation points.
It is in the spirit of Theorem~5.1 in~\cite{fuvi77}, where $f$ is
assumed to be locally compact and quasibounded at $0$.
Notice that the Leray--Schauder degree cannot be used here,
since we do not assume $f$ to be locally compact.

\begin{theorem} \label{teo-bif}
Let $\lambda_0, \lambda_1 \in \K$.
Assume $\lambda_0 \in \sigma (f,0) \backslash \sigma_\pi (f,0)$ and
$\lambda_1 \not \in \sigma (f,0)$.
Then, $\sigma_\omega (f,0) \cup \B (f,0)$ separates $\lambda_0$ from
$\lambda_1$; that is, $\lambda_0$ and $\lambda_1$ belong to different
components of $\K \backslash (\sigma_\omega (f,0) \cup \B (f,0))$.
\end{theorem}

\begin{proof}
Let $\lambda: [0,1] \to \K$ be a continuous path with
$\lambda (0) = \lambda_0$ and $\lambda (1) = \lambda_1$.
We need to prove that
$\lambda(t) \in \sigma_\omega (f,0) \cup \B(f,0)$
for some $t\in(0,1)$.

Assume by contradiction that
$\lambda(t) \not \in \sigma_\omega (f,0) \cup \B(f,0)$ for all
$t \in [0,1]$.
Then, there exists $r>0$ such that the equation
\[
\lambda(t)x-f(x), \quad t \in [0,1], \|x\| \leq r
\]
has only the trivial solutions $(\lambda(t),0)$, and we have
$\omega ((\lambda(t)-f)|_{\overline{V}}) >0$
for all $t \in [0,1]$, where $V= B_r$.
Consequently, the map $\lambda(t)-f$ is $0$-admissible on $V$ for all
$t$.
Since $\lambda_1 \not \in \sigma (f,0)$, we may assume that $V$ is so
small that $\lambda_1-f$ is $0$-epi on $V$.
Moreover, since $\lambda_0 \in \sigma (f,0)$,
$\lambda_0-f$ is not $0$-epi on $V$.
Thus, the two sets
\[
A_0= \{t \in [0,1] : \lambda(t)-f \mbox{ is not $0$-epi on $V$} \}
\quad \mbox{ and } \quad
A_1= \{t \in [0,1] : \lambda(t)-f \mbox{ is $0$-epi on $V$} \}
\]
are nonempty.
Clearly, $[0,1]$ is the disjoint union of $A_0$ and $A_1$.
On the other hand, as a consequence of Theorem~\ref{teo-stabilita},
$A_0$ and $A_1$ are open in $[0,1]$.
This is a contradiction because $[0,1]$ is connected, and the
assertion follows.
\end{proof}

\begin{corollary} \label{cor-bif}
Let $\lambda_0 \in \sigma (f,0) \backslash \sigma_\pi (f,0)$, and
assume that $\sigma (f,0)$ is bounded.
Then, the connected component of
$\K \backslash (\sigma_\omega (f,0) \cup \B (f,0))$
containing $\lambda_0$ is bounded.
\end{corollary}

\medskip

The next result is a sharpening of
Proposition~\ref{prop-componente-sigma}.

\begin{proposition} \label{prop-componente-bif}
Let $U$ be an open subset of $E$, with $\dim E = + \infty$.
Let $f \in C(U,E)$ be locally compact and assume that $\sigma (f,0)$
is bounded.
Then, $0 \not \in \Sigma (f,0)$ implies that the connected component
of $\K \backslash \B (f,0)$ containing $0$ is bounded.
\end{proposition}

\begin{proof}
Recall that $\sigma_\pi (f,0) = \{0\} \cup \Sigma (f,0)$
since $f$ is locally compact
(see Proposition~\ref{prop-spettro-compatte})
and that $\Sigma (f,0)$ is closed.
Thus, by Proposition~\ref{prop-componente-sigma}, if $|\lambda| >0$ is
sufficiently small then
$\lambda \in \sigma (f,0) \backslash \sigma_\pi (f,0)$.
Moreover, $\lambda \not \in \sigma (f,0)$ when $|\lambda|$ is
sufficiently large.
Hence, the assertion follows from Theorem~\ref{teo-bif}.
\end{proof}


\section{Illustrating examples}
\label{sect-examples}
\setcounter{equation}{0}

In this section we give some examples illustrating our main results.

We consider first the case $E = \C$.
Since $\C$ is finite dimensional, given $f: \C \to \C$ and
$p \in \C$, we have $\sigma_\pi (f,p) = \Sigma (f,p)$.
If, in addition, $f$ is positively homogeneous, then
\[
\Sigma(f,0) =
\{ \lambda \in \C : \lambda z = f(z) \mbox{ for some } z \neq 0 \}
= \B(f,0)
\]
(see Remark~\ref{oss-autovalori}).

\begin{example}
Let $f:\C \to \C$ be defined as $f(x+iy) = |x|+iy$.
The map $f$ is positively homogeneous and, consequently,
$\lambda = a+ib$ belongs to $\Sigma (f,0)$ if and only if the equation
\[
(a+ib)(x+iy)-(|x|+iy) =0
\]
admits a solution $x+iy$ in $S\sp 1$; that is, if and only if the
system
\[
\left\{
\begin{array}{l}
ax - |x| -by =0 \\
bx+ (a-1)y=0
\end{array}
\right.
\]
admits a nontrivial solution $(x,y)$.
An easy computation shows that $\Sigma (f,0) = S\sp 1$.
Observe that $d_0(f) = |f|_0 =1$, and this implies that the spectrum
$\sigma (f,0)$ is bounded.
Moreover, $f$ is not zero-epi at $0$.
To see this notice that, given $w \in \C$ with negative real part,
the equation $f(z) = w$ has no solutions.
Finally, since $\partial \sigma (f,0) \subseteq \Sigma (f,0)$,
we conclude that
$\sigma (f,0) = \{ a+ib : a\sp 2+b\sp 2 \leq 1 \}$
and $\Sigma (f,0) =S\sp 1$.
\end{example}

\begin{example}
Let $f: \C \to \C$ be defined as $f(x+iy) = sx+ty+i(ux+vy)$,
where $s,t,u,v$ are given real constants.
The map $f$ is positively homogeneous,
and linear if regarded from $\R \sp 2$ into itself.
Consequently, $\lambda = a+ib$ belongs to $\Sigma (f,0)$ if and only
if the equation
\[
(a+ib)(x+iy)-(sx+ty+i(ux+vy)) =0
\]
admits a solution $x+iy$ with $x\sp 2+y\sp 2 >0$; that is, if and only if
the homogeneous linear system
\[
\left\{
\begin{array}{l}
(a-s)x - (b+t)y =0 \\
(b-u)x + (a-v)y=0
\end{array}
\right.
\]
admits a nontrivial solution $(x,y)$.
This fact is equivalent to the condition
\[
\det \,
\left(
\begin{array}{cc}
a-s & -(b+t) \\
b-u & a-v
\end{array}
\right)
= 0,
\]
from which we get
\[
a\sp 2+b\sp 2-(s+v)a-(u-t)b+sv-tu=0.
\]
This is the equation of the circle $S_0$, centered at
$\left( \frac{s+v}2 , \frac{u-t}2 \right)$
with radius
\[
r= \sqrt{\frac{(s+v)\sp 2}4 + \frac{(u-t)\sp 2}4 -sv+tu}.
\]
Observe that $\sigma (f,0) = \Sigma (f,0)$.
Indeed, assume that $\lambda = a+ib$ does not belong to
$\Sigma (f,0)$, that is
\[
\det \,
\left(
\begin{array}{cc}
a-s & -(b+t) \\
b-u & a-v
\end{array}
\right)
\neq 0.
\]
This implies that $\lambda -f$ is a linear isomorphism as a map from
$\R \sp 2$ into itself.
In particular, $\lambda -f$ is a homeomorphism as a map from $\C$ into
$\C$.
It follows $\lambda \not \in \sigma (f,0)$ in view of
Proposition~\ref{prop-omeoloc}.
Hence, the whole spectrum $\sigma (f,0)$ coincides with the circle
$S_0$.

Notice that the spectrum reduces to a point (i.e.\ $r=0$) if and only
if $s=v$ and $t=-u$; that is, if and only if $f$ is linear as a
complex map.
\end{example}

\begin{example}
Let $g: \C \to \C$ be defined by
$g(x+iy) = \sqrt{x\sp 2+y\sp 2} + iy\sp n$, with $n \geq 2$.
As a consequence
of Theorem~\ref{teo-decomp-spettrale}-(3), we have
$\sigma (g,0) \equiv \sigma (f,0)$, where $f: \C \to \C$ is the
positively homogeneous map defined as
$f(x+iy) = \sqrt{x\sp 2+y\sp 2}$.
Now, it is not difficult to prove that
$\Sigma (f,0) = S \sp 1$.
Moreover, $\sigma (f,0) = \{ a+ib : a\sp 2+b\sp 2 \leq 1 \}$
since $0$ is not in the interior of the image of $f$
(so $f$ is not zero-epi at $0$).
Consequently, $\Sigma (g,0) = S\sp 1$ and
$\sigma (g,0) = \{ a+ib : a\sp 2+b\sp 2 \leq 1 \}$.
Hence, by Theorem~\ref{teo-bif},
we get $\B (g,0) = \Sigma (g,0) = S\sp 1$.
\end{example}

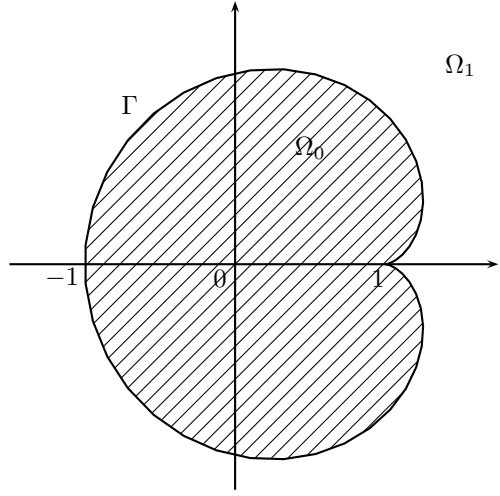
\begin{figure}[!tbhp]
\begin{center}
\begin{pspicture}(-3,-3)(3.5,3.5)
\psaxes[ticks=none,labels=none]{->}(0,0)(-3,-3)(3.5,3.5)
\parametricplot[fillstyle=hlines,hatchcolor=black,hatchwidth=.4pt]{0}{360}%
{1 t cos t cos mul sub t cos add 2 mul %
t sin t cos mul t sin sub 2 mul}
\rput[b](1.0,1.4){$\Omega_0$}
\rput[b](3.0,2.5){$\Omega_1$}
\rput[b](-1.4,2){$\Gamma$}
\rput[b](-2.3,-0.3){$-1$}
\rput[b](1.9,-0.3){$1$}
\rput[b](-0.2,-0.3){$0$}
\end{pspicture}
\end{center}
\caption{The spectrum of $f: \C \to \C$,
$x+iy \mapsto \sqrt{x\sp 2+y\sp 2}+iy$.}
\label{fig-cardio}
\end{figure}

\begin{example} \label{es-cardio}
Let $f:\C \to \C$ be defined as
$f(x+iy) = \sqrt{x\sp 2+y\sp 2}+iy$.
Since $f$ is positively homogeneous,
one can show that $a+ib \in \Sigma (f,0)$ if and only if
\[
(a-1)\sp 2 + b\sp 2 = (a\sp 2 + b\sp 2 -a)\sp 2,
\]
which is the equation of a closed curve $\Gamma$ (a cardioid).
The curve $\Gamma$ divides the complex plane in two connected
components, $\Omega_0$ (containing $0$) and
$\Omega_1$ (unbounded).
Clearly, $\lambda \not \in \sigma (f,0)$ if $\lambda$ belongs to
$\Omega_1$.
Furthermore, $\lambda \in \sigma (f,0)$ for any $\lambda \in \Omega_0$
since $f$ is not zero-epi at $0$.
Hence,
\[
\sigma (f,0) = \overline \Omega_0 = \Omega_0 \cup \Gamma
\mbox{ and } \Sigma (f,0) =\Gamma
\]
(see Figure~\ref{fig-cardio}).
\end{example}

\begin{figure}[!htbp]
\begin{center}
\begin{pspicture}(-3,-3)(4,4)
\pscircle[fillstyle=hlines,hatchcolor=black,hatchwidth=.4pt](0.75,0){2.25}
\pscircle[fillstyle=solid,hatchcolor=white,hatchwidth=.4pt](2.25,0){0.75}
\pscircle[fillstyle=hlines,hatchcolor=white,hatchwidth=.4pt](2.25,0){0.75}
\rput[b](3.5,2.1){$\Omega_2$}
\rput[b](2.2,0.1){$\Omega_1$}
\rput[b](0.8,1.2){$\Omega_0$}
\rput[b](-1.2,1.5){$S_+$}
\rput[b](1.8,-1.0){$S_-$}
\rput[b](-0.2,-0.3){$0$}
\rput[b](-1.8,-0.5){$-\frac12$}
\rput[b](3.1,-0.3){$1$}
\rput[b](1.4,-0.5){$\frac12$}
\psaxes[ticks=none,labels=none]{->}(0,0)(-2.5,-3)(4,3.5)
\end{pspicture}
\end{center}
\caption{The spectrum of $f: \C \to \C$,
$x+iy \mapsto \frac{|x|}2+iy$.}
\label{fig-tang}
\end{figure}
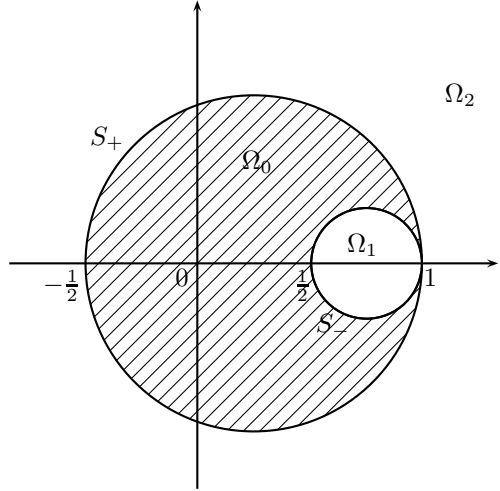

\begin{example} \label{es-tang}
Let $f: \C \to \C$ be defined as $f(x+iy) = \frac{|x|}2+iy$.
Notice that $d_0(f) = \frac12$ and $|f|_0 =1$.
Since $f$ is positively homogeneous,
one can show (see e.g.~\cite{fuvi75})
that $\Sigma (f,0)$ is the union of two circles:
$S_+ = \left\{ \lambda \in \C :
\left| \lambda-\frac14 \right| = \frac34 \right\}$
and
$S_- = \left\{ \lambda \in \C :
\left| \lambda-\frac34 \right| = \frac14 \right\}.$
Consequently, $\C \backslash (S_+ \cup S_-)$ consists of three
connected components, $\Omega_0$ (containing $0$),
$\Omega_1$ (surrounded by $S_-$) and
$\Omega_2$ (unbounded).
One can check that, when $\lambda$ belongs to
$\Omega_1 \cup \Omega_2$, the map $\lambda -f$ is a local
homeomorphism around zero.
Thus, $\lambda \not \in \sigma (f,0)$ on the basis of
Proposition~\ref{prop-omeoloc}.
Moreover, since $f$ is not zero-epi at $0$, we get
$\lambda \in \sigma (f,0)$ for all $\lambda \in \Omega_0$.
Hence,
\[
\sigma (f,0) = \overline \Omega_0 = \Omega_0 \cup (S_+ \cup S_-)
\mbox{ and }
\Sigma (f,0) =S_+ \cup S_-
\]
(see Figure~\ref{fig-tang}).
\end{example}

\medskip

We close with two examples in the infinite dimensional context.
The first one regards a differentiable interesting map,
and the second one deals with a merely continuous map.

\begin{example}
Let $E$ be a real Hilbert space
(of dimension greater than $1$), and consider the nonlinear map
$f(x) = \| x \| x$.
Observe that $f$ is Fr\'echet differentiable at any $p \in E$ with
\[
f'(p) v = \| p \| v + \frac{(p,v)}{\|p\|} \, p,
\quad v \in E
\]
if $p \neq 0$, and $f'(0) =0$.
Hence, by Corollary~\ref{cor-equivalenza} we have
$\sigma (f,p) = \sigma (f'(p))$.
Given $p \neq 0$, in order to compute $\sigma (f'(p))$, observe that
$f'(p)$ is of the form
$L = c I + K$, with $c \in \R$ and $K : E \to E$
with finite dimensional image, say $E_0$
(in this case $\dim E_0 = 1$).
Since $\sigma (L) = \{ c \} \cup \sigma (L_0)$,
where $L_0$ denotes the restriction of $L$ to $E_0$,
we get
$\sigma (f,p) = \sigma (f'(p)) = \{ \|p\| \} \cup \{ 2\|p\| \}$
if $p \neq 0$ and, clearly,
$\sigma (f,0) = \sigma (f'(0)) = \{ 0 \}$.
\end{example}

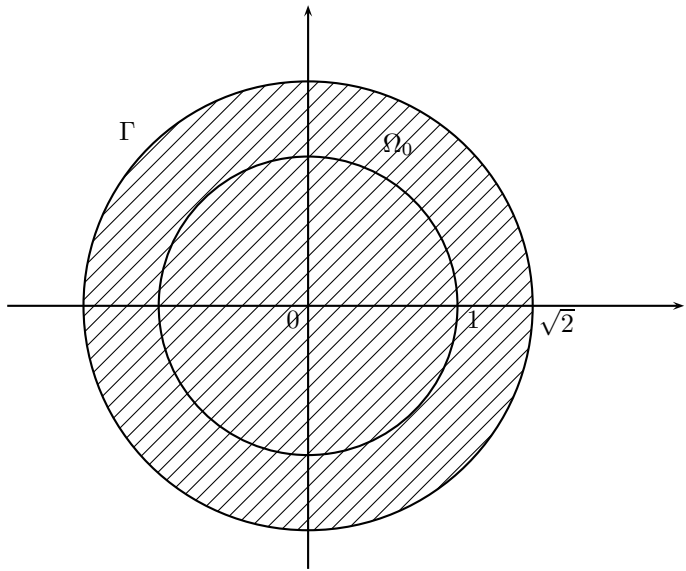
\begin{figure}[!htbp]
\begin{center}
\begin{pspicture}(-3.5,-3.5)(4,4)
\psaxes[ticks=none,labels=none]{->}(0,0)(-4,-3.5)(5,4)
\pscircle(0,0){2}
\pscircle[fillstyle=hlines,hatchcolor=black,hatchwidth=.4pt](0,0){3}
\rput[b](1.2,2,2){$\Omega_0$}
\rput[b](-2.4,2.2){$\Gamma$}
\rput[b](3.3,-0.4){$\sqrt2$}
\rput[b](2.2,-0.3){$1$}
\rput[b](-0.2,-0.3){$0$}
\end{pspicture}
\end{center}
\caption{The spectrum of $f: \ell \sp 2 (\C) \to \ell \sp 2 (\C)$,
$z \mapsto (\|z\|,z_1,z_2,z_3,\dots)$.}
\label{fig-concent}
\end{figure}

\begin{example} \label{es-concent}
Let $f: \ell \sp 2 (\C) \to \ell \sp 2 (\C)$ be defined by
\[
f(z)=(\|z\|,z_1,z_2,z_3,\dots),
\]
where $z=(z_1,z_2,z_3,\dots)$.
Notice that $f$ is positively homogeneous, and is the sum of the
\textit{right-shift operator}
$L: \ell \sp 2 (\C) \to \ell \sp 2 (\C)$, defined as
$Lz=(0,z_1,z_2,z_3,\dots)$,
and the finite dimensional map
$k: \ell \sp 2 (\C) \to \ell \sp 2 (\C)$, defined as
$k(z)=(\|z\|,0,0,0,\dots)$.

Let us compute $\sigma (f,0)$.
An easy computation shows that $d(f) = |f| = \sqrt 2$.
Moreover, $\alpha(f) = \omega(f) =1$.
Indeed, since $k$ is compact and $f=L+k$, we have
$\alpha(f) = \alpha(L)$ and $\omega(f) = \omega(L)$.
Now, $\alpha(L) = \omega(L) =1$, $L$ being an isometry between the
space $\ell \sp 2 (\C)$ and a subspace of codimension one.
Therefore, Proposition~\ref{prop-stime} implies
$\sigma_\omega (f,0) \subseteq \{ \lambda \in \C: |\lambda| = 1 \}$
and
$\Sigma (f,0) \subseteq \{ \lambda \in \C: |\lambda| = \sqrt2 \}$.
Let us show that the converse inclusions hold.

First, let us prove that
$\sigma_\omega (f,0) = S\sp 1$.
Since $\omega (\lambda - f) = \omega (\lambda - L)$, it is enough to
show that $\omega (\lambda - L) =0$ when $|\lambda| =1$.
To this end, recall that a linear operator $T$ is left semi-Fredholm
if and only if $\omega(T) >0$.
Thus, $\lambda - L$ is left semi-Fredholm for $|\lambda| \neq 1$.
Recall also that the index of $\lambda - L$,
\[
\ind (\lambda - L) =
\dim \Ker (\lambda - L) - \dim \coker (\lambda - L)
\in \{ - \infty \} \cup \Z,
\]
depends continuously on $\lambda$. Therefore, it is constant on any
connected set contained in $\C \backslash S\sp 1$.
This implies that $\ind (\lambda - L) =-1$ when $|\lambda| <1$ since
$\ind(-L) =-1$.
On the other hand, $\ind(\lambda - L) =0$ if $|\lambda| >1$ since, as
well known, $\sigma(L) = \{\lambda \in \C: |\lambda| \leq 1 \}$.
Thus, the subset $\sigma_\omega (L)$ of $S\sp 1$ separates the two
open sets
$\{\lambda \in \C: |\lambda| <1 \}$
and
$\{\lambda \in \C: |\lambda| >1 \}$.
Consequently, $\sigma_\omega (L) = S\sp 1$.
Hence, $\sigma_\omega (f,0) = S\sp 1$.

To show that
$\Sigma (f,0) = \{ \lambda \in \C: |\lambda| = \sqrt2 \}$,
assume $|\lambda| =\sqrt2$.
Since $\lambda \not \in \sigma_\omega (f,0)$,
Proposition~\ref{prop-autovalori} implies that
$\lambda \in \Sigma (f,0)$ if and only if
$\lambda$ is an eigenvalue of $f$ at $0$, i.e.,
\[
\min_{\|z\|=1} \|\lambda z - f(z)\|=0.
\]
Simple computations show that this condition is satisfied
when $|\lambda| = \sqrt2$.

As a consequence of the above arguments, $\sigma_\pi (f,0)$ is the
union of two circles centered at the origin.
Now, observe that
$q_0 (f) = \max\{ \alpha_0 (f),|f|_0\} = \sqrt2$ and hence
$\sigma (f,0) \subseteq \{ \lambda \in \C: |\lambda| \leq \sqrt2 \}$.
Let us prove that
$\sigma (f,0) = \{ \lambda \in \C: |\lambda| \leq \sqrt2 \}$.
For this purpose remind that, if $W$ is a connected component of
$\C \backslash \sigma_\pi(f,0)$, then the maps of the form
$\lambda-f$, with $\lambda \in W$, are either all zero-epi at $0$ or
all not zero-epi at $0$ (see Corollary~\ref{cor-componenti-sigma}).

We claim that,
if $|\lambda| <1$, then $\lambda - f$ is not zero-epi at $0$.
Set $e_1 = (1,0,0, \dots)$.
Observe that, given $\eps >0$, the equation $f(z) = -\eps e_1$ has no
solutions. Consequently, $f$ is not zero-epi at $0$ and the claim is
proved.

Let us show that $\lambda - f$ is not $0$-epi at $0$
also for $1< |\lambda| <\sqrt2$.
Indeed, fix $\lambda \in \C$ with $1< |\lambda| <\sqrt2$.
We claim that, given $\eps >0$, the equation
\begin{equation} \label{equ-es}
\lambda z - f(z) = \eps e_1, \quad z \in E
\end{equation}
has no solutions.
Recall that $\lambda - L$ is an isomorphism and set
$v_{\lambda} = (\lambda - L) \sp{-1} (e_1)$.
Since $f$ is the sum of $L$ and $k$,
and the image of $k$ lies in the subspace spanned by $e_1$,
the solutions of~\eqref{equ-es} lie in the one
dimensional subspace $E_{\lambda}$ spanned by $v_{\lambda}$.
Therefore, any solution of~\eqref{equ-es} is of the type
$\xi v_{\lambda}$, $\xi \in \C$.
An easy computation shows that
$\|v_{\lambda}\|\sp 2 = \frac1{|\lambda|\sp 2-1}$.
Thus,
\[
\lambda (\xi v_{\lambda}) - f(\xi v_{\lambda}) =
\xi e_1 - \|\xi v_{\lambda}\| e_1 =
\left( \xi - |\xi| \frac1{\sqrt{|\lambda|\sp 2-1}} \right) e_1.
\]
Consider now the equation
\begin{equation} \label{equ-es2}
\xi - |\xi| \frac1{\sqrt{|\lambda|\sp 2-1}} = \eps,
\quad \xi \in \C
\end{equation}
which is equivalent to~\eqref{equ-es}.
It is not difficult to see that equation~\eqref{equ-es2} has no
solutions when $1< |\lambda| <\sqrt2$.
Consequently, equation~\eqref{equ-es} has no solution, as claimed.
Hence, $\lambda - f$ is not zero-epi at $0$
when $1< |\lambda| <\sqrt2$.

{}From the above discussion we get
\[
\sigma(f,0) = \{\lambda \in \C: |\lambda| \leq \sqrt2 \}.
\]

\smallskip

Consider now any map $h: \ell \sp 2 (\C) \to \ell \sp 2 (\C)$ which is
compact and such that $h(z) = o (\| z\|)$ as $\|z\| \to 0$, and let
$g=f+h$. Then, as a consequence of
Theorem~\ref{teo-decomp-spettrale}-(3), we have
$\sigma (g,0) \equiv \sigma (f,0)$.
In particular,
\[
\sigma(g,0) = \{\lambda \in \C: |\lambda| \leq \sqrt2 \}=\Omega_0.
\]
Moreover,
\[
\sigma_\omega (g,0) = S\sp 1
\quad \mbox{ and } \quad
\Sigma(g,0) = \{\lambda \in \C: |\lambda| = \sqrt2 \} =\Gamma
\]
(see Figure~\ref{fig-concent}).
Hence, from Theorem~\ref{teo-bif}, it follows that any
$\lambda \in \C$ with $|\lambda|=\sqrt2$ is a bifurcation point for
$g$. That is,
\[
\B(g,0) = \{\lambda \in \C: |\lambda| = \sqrt2 \} =\Gamma.
\]
In view of Theorem~\ref{teo-bif}, this stability of the set of
bifurcation points depends on the fact that $\Sigma(f,0)$ locally
separates $\sigma(f,0)$ from its complement.

Notice that in the above example we have detected a bifurcation
phenomenon that cannot be investigated via the classical
Leray--Schauder degree theory.
Moreover, since the map $f$ is a compact perturbation of a linear
Fredholm operator of negative index, also the more recent
degree theory for compact perturbations of Fredholm
operators of index zero (see~\cite{befu} and references therein)
cannot be applied.
\end{example}

\bigskip


\end{document}